%% file: paper_arxiv.tex
\documentclass[11pt,leqno]{article}

\usepackage{natbib}
\usepackage{graphicx}
\usepackage{a4wide}
\usepackage{multirow}
\usepackage{bm}
\usepackage{nicefrac}
\usepackage{algorithm,algorithmic}
\usepackage{csquotes}
\usepackage{soul}
\usepackage[inkscapelatex=true]{svg}
\usepackage{amsmath}
\usepackage{amssymb}

\usepackage{tikz}
\usepackage{tikz-3dplot}
\tdplotsetmaincoords{70}{30}
\usetikzlibrary{shapes.geometric, arrows, fit, positioning}
\usepackage{subcaption, graphicx}
\tikzstyle{startstop} = [rectangle, rounded corners, minimum width=1cm, minimum height=0.3cm, text centered, draw=black, fill=red!50]
\tikzstyle{process} = [rectangle, minimum width=1cm, minimum height=0.3cm, text centered, draw=black, fill=orange!50]
\tikzstyle{morprocess} = [rectangle, minimum width=1cm, minimum height=0.3cm, text centered, draw=black, fill=green!50]
\tikzstyle{decision} = [diamond, aspect=4, minimum height=0.1cm, text centered, inner sep=0.1cm, draw=black, fill=blue!50]
\tikzstyle{decision2} = [diamond, aspect=3, minimum height=0.1cm, text centered, inner sep=0.001cm, draw=black, fill=blue!50]
\tikzstyle{arrow} = [thick,->,>=stealth]

\newcommand{\enthalpy}{\ensuremath{T}}
\newcommand{\adjointVar}{\ensuremath{T_{a}}}

\newcommand{\real}{\mathbb{R}}
\newcommand{\<}{\leq}
\newcommand{\e}{{\bm e}}
\newcommand{\fomtau}{\ensuremath{\tau_{fom}}}
\newcommand{\romtau}{\ensuremath{\tau_{mor}}}
\newcommand{\krylovCommand}[3]{\ensuremath{\mathcal{K}_{#3}\left(#1;\ #2\right)}}
\newcommand{\norm}[1]{\ensuremath{\left\|#1\right\|}}
\newcommand{\state}[2]{\ensuremath{#1^{\left(#2\right)}}}
\newcommand{\cgBasis}{\ensuremath{W}}

\DeclareMathOperator*{\argmin}{arg\,min}

\begin{document}

\title{Comparison of model order reduction techniques with one-shot procedure for topology optimization for thermal applications
}
    \author{Luis Fernando Cusicanqui Lopez\footnote{Corresponding author. \\ email: luisfernando.cusicanquilopez@kuleuven.be},
        Ramadan Krasniqi, 
        Florian Feppon, 
        Karl Meerbergen
}
\maketitle

\begin{abstract}
Density-based topology optimization has become a powerful method for automatically generating optimized designs in
    a wide variety of applications.
    However, it comes with a large computational cost when 
    solving the physical model requires large-scale simulations.  Here, we investigate the use of model order reduction (MOR) techniques to accelerate the simulations in the context of thermal design applications. 
    We project the governing and the adjoint equations onto a low-dimensional subspace by constructing two distinct reduced bases---one for the forward state and one for the adjoint system---using solution snapshots from previous design iterations. These snapshots are generated using either the high-fidelity solver or inaccurate fast solvers, such as the one-shot method \citep{amir2024one}. Additionally, we demonstrate that properly selecting the stopping criterion for the iterative linear solver is crucial for the effective use of reduced models.
    In our 3D example, the proposed framework reduces the overall total simulation time relative to the high-fidelity workflow by a factor up to $3$ when combined with high-fidelity solves and a factor up to $16$ when combined with the one-shot method. Moreover, we find that the reduced order model approach is able to achieve a speed up of  $1.54$ with respect to the one-shot method. 
\end{abstract}

\section{Introduction}
\label{s:introduction}
Topology Optimization (TO) is a powerful numerical method that enables the automatic generation of optimized designs for a given objective \citep{bendsoe2013topology}. 
It is an iterative design process that involves solving equations that model the underlying physics at each iteration. However, this task is computationally expensive for three-dimensional, large-scale applications, as it requires solving  large linear systems \citep{amir2014multigrid,MARTINEZFRUTOS201747,aage2015topology,amir2010efficient}. Consequently, the solution process often becomes the bottleneck of the overall design cycle. Therefore, reducing the computational cost associated with simulating these large linear systems is critical to speeding up the entire design process. \\
Several approaches have been explored to address this bottleneck. These include high performance parallel computing \citep{aage2015topology,aage2017giga}, reducing the number of design variables \citep{guest2010reducing} and reducing the number of solver iterations by using an efficient preconditioner as in \citep{amir2014multigrid} or applying subspace recycling techniques \citep{choi2019acceleratingdesignoptimizationusing,parks2006recycling,wang2007large}. \\
% In this paper, we focus specifically on reducing the computational cost associated with simulations at each iteration of the optimization process. \\
The cost associated with the simulations is dominated by the cost of solving the resulting large linear systems during the design process. These linear systems are obtained by discretizing the governing equations, where the design variable is typically embedded in the material properties. As the optimization progresses, the change in these design variables becomes smaller, leading to small differences between consecutive linear systems. This can be exploited by numerical methods to reduce the computational costs incurred by solving large-scale linear systems. A promising approach is projection-based model order reduction that builds a reduced basis to project the large-scale linear system onto a low dimensional system and generate approximate solutions. This basis is typically built using solutions corresponding to previous optimization steps during the design process. In topology optimization, this approach has been initially followed by \citep{GoguMorTopopt,xiao:hal-02531425} that used Gram-Schmidt and POD to build the basis and applied model order reduction to the forward equation for a minimal compliance structural problem, however the methodology is limited to self adjoint problems.
\citep{choi2019acceleratingdesignoptimizationusing} developed a model order reduction framework that can also be applied for non self-adjoint optimization problems. However, the incremental algorithm used to update the bases fails to correctly recover the desired subspace once the maximum number of basis vectors is reached. This incremental approach orthogonalizes new snapshots with respect to the present basis vectors from previous iterations. Once the maximum basis size is reached, the algorithm removes the oldest orthogonalized basis vector and orthogonalizes the new incoming snapshot against the remaining basis vectors. Subsequent basis vectors were originally orthogonalized against this discarded vector; therefore, information from the original snapshots that was present in the removed vector is permanently lost. Consequently, the remaining basis can no longer accurately reconstruct the desired subspace, i.e., the original snapshots. \citep{QIAN2022110917} applied a dual reduction, where the forward and the adjoint equations were reduced separately with a different basis, to speedup the design process for a time-dependent structural optimization problem. 
An alternative strategy to mitigate computational costs involves solving the state equations inaccurately to reduce the number of solver iterations needed at each optimization step. For instance, \citep{amir2010efficient} solves the state equations inaccurately by dynamically adjusting the residual threshold of the iterative solver over the course of the optimization. More recently, in \citep{amir2024one}, the authors used a one-shot procedure for a minimum compliance structural problem. In this context, one-shot means a single iteration of the conjugate gradient method.
This resulted in large speedups without jeopardizing the quality of the final design. \\
There is a conflict between computational speed and simulation accuracy achieved by these methods. This trade-off is governed by the linear solver's stopping criterion. If the criterion is too soft, the result may be inaccurate and lead to bad results; if it is too stringent, speedup is lost. Furthermore, we identify that the standard choice of stopping criterion for conjugate gradient (CG) is too strict when applied to the adjoint equation leading to a large number of solver iterations or an inefficient use of the reduced model. To the best of our knowledge, alternative ways of computing the residual norms have not been studied before in the context of topology optimization. Therefore, based on results from \citep{banoczi1998lack,arioli1992stopping}, we propose an alternative residual norm as stopping criterion that leads to a more effective use of the reduced models. \medskip \\
In this paper, we investigate the use of MOR to mitigate the computational cost of simulations in the context of topology optimization for heat conduction problems.
More specifically, we reduce the forward and the adjoint equation separately with their own corresponding bases, which are constructed using converged solutions from high-fidelity solvers and inaccurate solutions from the one-shot method.
First, we propose a Gram-Schmidt based incremental algorithm to correctly recover the desired subspace when removing `old' vectors from the basis once its maximum dimension is reached.
Second, we propose an alternative stopping criterion from \citep{banoczi1998lack,arioli1992stopping} which leads to a more uniform use of the forward and adjoint reduced models improving the overall performance. Finally, we notice that the computational cost of a full multigrid precondition step at every design iteration is still higher than one model order reduction (MOR) step. Therefore, we investigate whether MOR bases constructed from these inaccurate one-shot solutions are robust enough to achieve a qualitative design. This combination allows us to exploit the low computational cost of MOR and the efficiency of the one-shot method. The one-shot method \citep{amir2024one} already gives a speedup of $10$ for the problem considered. In combination with MOR, an additional speedup of $1.53$ can be achieved. \medskip \\
This article is organized as follows: Section \ref{ss:topopt} introduces the
necessary background on topology optimization for a classical heatsink design problem. Section \ref{s:linearSolvers} describes the choice of the linear solver and the introduction of the new stopping criterion used, followed by the proposed MOR framework in \ref{s:mor}. Finally, Section \ref{s:experiments}  shows two- and three-dimensional numerical results using the proposed workflow, and provides a comparison with the available state-of-the-art techniques. 

\section{Density-based topology optimization}\label{ss:topopt}
Topology optimization (TO) for thermal applications seeks to find an optimal material distribution minimizing an objective function subject to specific constraints \citep{bendsoe2013topology}. In the density based TO framework, the material distribution in the design domain is represented using a discrete variable $\rho\left(x\right)$, taking values $0$ or $1$, where $0$ corresponds to the original material and $1$ to the design material. Since this formulation leads to a large-scale integer optimization problem, it is customarily relaxed by allowing intermediate (continuous) values for the design variable $\rho$, converting the design problem into a continuous optimization problem. A density-based topology optimization problem involving a simple volume constraint can be formulated as a constrained optimization program:
\begin{equation}
	\begin{split}
		\min_\rho \mathcal{J}(\rho, \enthalpy), \\
		\mathrm{s.t.} \\
	   \left\{\begin{split}
	       R(\rho, \enthalpy) = 0,\\
		\frac{\int\rho d\Omega}{\int d\Omega} \le v_{\mbox{frac}}, \\
		0 \le \rho \le 1,  
	   \end{split}	\right . 
	\end{split}
\end{equation}
where $\mathcal{J}(\rho, \enthalpy)$ is the objective function and $R(\rho, \enthalpy)$ represents the governing equations that the design must satisfy. The maximum allowed relative volume of the design is imposed to be $v_{\mbox{frac}}$. In this paper, we consider thermal conduction within a solid material where the physical behavior is governed by the steady-state heat conduction equation for the temperature field $\enthalpy$ in the domain  given by
\begin{equation}\label{eq:energyEqn}
    \left\{\begin{split}
        -\nabla\cdot\kappa(\rho)\nabla \enthalpy &= Q, \text{ in }\Omega \\
        \enthalpy &=  \enthalpy_d \;\; \text{in} \;\;  \partial \Omega_d, \\
        \nabla \enthalpy &= q_n \;\; \text{in} \;\; \partial \Omega_n, 
    \end{split}\right.
\end{equation}
where $\Omega$ is the design domain, $\kappa(\rho)$ the heat conductivity and $Q$ is an external heat source. The boundary $\partial\Omega$ is split into a Dirichlet boundary $\partial \Omega_d$ where we a fixed temperature $\enthalpy_d$ is imposed, and a Neumann boundary  $\partial \Omega_n$ where the heat flux is specified $q_n$. The heat conductivity $\kappa(\rho)$ is determined by the design variable $\rho$ following the solid isotropic material with penalization (SIMP) rule as in \citep{gersborg2006topologyHeat}:
\begin{equation}\label{eq:SIMP}
	\kappa(\rho) = \kappa_{\min} + \rho^P(\kappa_{\max} - \kappa_{\min}),
\end{equation}
where $\kappa_{\min}$ and $\kappa_{\max}$ are the heat conductivities of the two
materials and $P$ is a penalization exponent that 
penalizes intermediate (gray) design values.\\
The use of any gradient-based optimization solver requires the computation of the
sensitivities of the objective and constraint functions with respect to the design
variable $\rho$. We adopt the optimize then discretize framework \citep{Hinze2009} where the sensitivities of the continuous equations are derived before being discretized for computing a decent direction. A standard derivation leads to the following expression for the Fréchet derivative of the objective with respect to the design variable:
\begin{equation}\label{eq:grad}
	\frac{\partial \mathcal{J}}{\partial \rho} \cdot \delta\rho =
   - \int_{\Omega}\frac{\partial \kappa}{\partial\rho} \nabla \enthalpy \cdot \nabla
    \adjointVar\,  \delta \rho \mathrm{d} \Omega,
\end{equation}
where $\adjointVar$ is the adjoint variable obtained by solving the adjoint equation, given in its weak form as:
\begin{equation}\label{eq:adjEnergyEqn}
	\int_\Omega \kappa(\rho)\nabla v \cdot \nabla \adjointVar\mathrm{d}\Omega  =
    \frac{\partial \mathcal{J}}{\partial T}\cdot v, \quad 
    \forall v\in \mathcal{H}^1_0(\Omega),
\end{equation}
where $\partial \mathcal{J}/\partial \rho$ and $\partial
\mathcal{J}/\partial T$ are the partial Fréchet derivatives of
$\mathcal{J}$, and where  $\mathcal{H}^1_0(\Omega)$ is the Sobolev space of square integrable functions
with weak derivatives in $L^2(\Omega)$ and vanishing on $\partial
\Omega_d$. \\
We apply a PDE regularizing filter to mitigate checkerboard  effects as described in \citep{pdelazarov2011filters,sigmund200199} where the density field is filtered by solving the differential equation
\begin{equation}\label{eq:pde_filter}
    -\lambda^2\nabla^2\tilde{\rho} + \tilde{\rho} =\rho,
\end{equation}
where $\lambda$ is the length parameter and $\tilde{\rho}$ is the filtered density. To determine the value of $\lambda$, we use the relation from \citep{pdelazarov2011filters}
\begin{equation}\label{eq:relation_filter}
    \lambda = \frac{1}{2\sqrt{3}}\lambda_h,
\end{equation}
where $\lambda_h$ is the radius of the support domain. \medskip \\
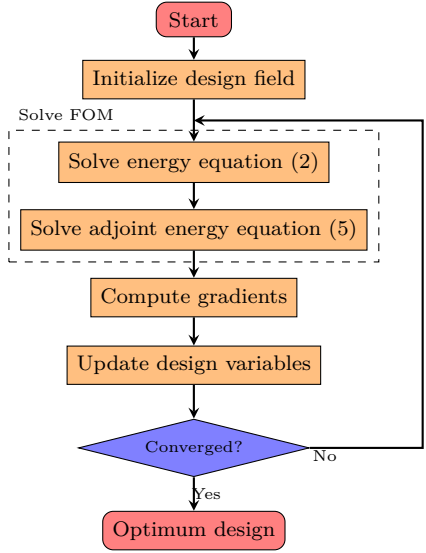
\begin{figure}[tp]
\begin{center}
    \input{workflowDiagram}
\end{center}
    \caption{Classic topology optimization workflow.}\label{fig:normalFlow}
\end{figure}
Figure~\ref{fig:normalFlow} presents a flow chart that illustrates all the steps
involved in a standard topology optimization procedure. First, the design
variable $\rho$ is assigned to an initial value. Then, every optimization iteration consists in the following three steps. First, the forward and adjoint governing
equations are solved, which together form the Full Order Model (FOM). This step is referred to as
\textbf{Solve FOM}. Second the gradient is computed using equation \eqref{eq:grad},
and third, the design is updated using an optimization algorithm.
In this paper, we use the \textit{Method of Moving Asymptotes} (MMA)
\citep{svanberg1987method} as the optimizer to update the design variables while
handling  the volume constraint in addition to
the bound constraints $0\<\rho\<1$. The
iterative procedure continues  after a fixed number of iterations set to $250$ in our implementation. In this procedure, the
most computationally expensive steps in each design iteration are the calls to the
linear system solvers of the governing forward and adjoint equations. 
\section{Linear system solvers}\label{s:linearSolvers}
The step \textbf{Solve FOM} requires to solve the forward and the adjoint equations. To this end, we have to discretize equations equations \eqref{eq:energyEqn} and \eqref{eq:adjEnergyEqn} at every design iteration. After discretization, two linear systems of the form
\begin{equation}\label{eq:linearSystem}
	A x = b
\end{equation}
need to be solved, where $A \in \real^{n \times n}$ is the coefficient matrix, $x \in \real^n$ is
a vector containing the field values of the mesh and $b \in \real^n$
is the discretized source term. The dimension $n$ is directly related to the
discretization resolution, which is typically large for topology optimization problems.
Since, equations \eqref{eq:energyEqn} and \eqref{eq:adjEnergyEqn} involve identical
operators, the coefficient matrices $A \in \real^{n \times
n}$ are identical for both, state and adjoint discretized equation, provided the same discretization methods are applied. In this paper, the conjugate gradient method is used to solve the governing equations.
\subsection{Preconditioned Conjugate Gradient method}
The Conjugate Gradient method (CG)
\citep{hestenes1952methods,saad2003iterative,gutknecht2007brief} is an iterative
method for solving the linear system $A x=b$ with  symmetric positive definite matrix
$A$. It generates an iterative $(x^{(j)})$ of
solutions that converge to the exact solution of the linear system, denoted by
$x^{*}$.
Let $r^{\left(0\right)}$ be the residual corresponding to the initial guess
$x^{(0)}$, defined as
\[
r^{\left(0\right)} = b - A x^{\left(0\right)} = A (x^* - x^{\left(0\right)}). 
\]
The residual $r^{(0)}$ generates  the Krylov subspace
\[
\krylovCommand{A}{r^{\left(0\right)}}{j} = \mbox{span}\left(r^{\left(0\right)}, Ar^{\left(0\right)}, \ldots, A^{j-1}r^{\left(0\right)}\right).
\]
Denote by $W\in\mathbb{R}^{n\times j}$ the matrix whose columns form a basis of
$\krylovCommand{A}{r^{\left(0\right)}}{j}$, namely $W\e_j=A^{j}r^{(0)}$.
The iterates 
$x^{(j)}$ generated by the CG method 
are of the form $\state{x}{j} = x^{\left(0\right)} + \widetilde{x}^{(j)}$, with
$\widetilde{x}^{(j)}$ the Galerkin projection with respect to the $A$-norm of
$x^*-x^{\left(0\right)}$ on the Krylov
subspace
\[
\widetilde{x}^{(j)} = \cgBasis \left(\cgBasis^TA\cgBasis \right)^{-1}\cgBasis^Tr^{\left(0\right)} \in \krylovCommand{A}{r^{\left(0\right)}}{j}.
\]
One of the main properties of this method is that the approximation
\state{x}{j} minimizes the $A$-norm, or energy norm of the error:
\begin{equation*}
	\widetilde{x}^{(j)} = \argmin_{\widetilde{x}^{(j)} \in\  \krylovCommand{A}{r^{\left(0\right)}}{j}}\norm{x - x^*}_{A} = \|A x-b\|_{A^{-1}},
\end{equation*}
where $\norm{x - x^*}_{A}^2= (x-x^*)^TA(x-x^*)$. The matrices arising from
high-resolution discretizations of physical models exhibit poor conditioning, which results in the need for preconditioners to achieve satisfactory convergence.
In this paper, we consider the geometric multigrid as preconditioner of the CG method. This is the standard for large-scale topology optimization problems involving the heat conduction or the linear elasticity model
\citep{amir2014multigrid,aage2017giga,aage2015topology,nobel20153d}. 
\subsection{Stopping criterion}\label{ss:stop}
The number of CG iterates is determined by a stopping criterion which is satisfied when a measure of the residual error is smaller than a certain tolerance. Additionally, we also use this error measure to assess the quality of our reduced order model: if this error is small, we do not need to solve the FOM, otherwise, the FOM is solved and the MOR is updated. Therefore, a proper stopping criterion is crucial for  global efficiency, to avoid unnecessary CG iterations if the criterion is too severe, while still not compromising accuracy.
Previous works investigating MOR in topology optimization relied on a stopping criterion based on the relative error of the residual \citep{choi2019acceleratingdesignoptimizationusing,GoguMorTopopt} defined as
\begin{equation}\label{eq:stop1}
w_1 = \frac{\|r^{(j)}\|}{\|b\|} \leq \tau ,
\end{equation}
where $\|r^{(j)}\|=\|b-Ax^{(j)}\|$ is the 2-norm of the residual after $j$ Krylov
iterations and $\tau$ is a user-defined absolute tolerance value. When this stopping criterion is chosen, the approximation error is bounded by \citep{arioli1992stopping}
\begin{equation}\label{eq:bound_stop_1}
    \frac{\mathrm{cond}(A,b)}{\mathrm{cond}(A)}w_1 \le \frac{\|x^*-x\|}{\|x^*\|}\le \mathrm{cond}(A,b)w_1 , 
\end{equation}
where $x^*$ is the exact solution, $\mathrm{cond}(A)$ is the condition number of matrix $A$ and $\mathrm{cond}(A,b)=\nicefrac{\|A^{-1}\|\|b\|}{\|x^*\|}$ is the condition number of $A$ that depends on the right hand side and fulfills the following condition
\begin{equation}
    1 \le \mathrm{cond}(A,b) \le \mathrm{cond}(A).
\end{equation}
An alternative stopping criterion is
\begin{equation}\label{eq:stop2}
w_2 = \frac{\|r^{(j)}\|}{\|b\| + \|A\| \|x^{(j)}\|} \leq \tau,
\end{equation}
where the relative error is now bounded by
\begin{equation}\label{eq:bound_stop_2}
    \left(1+\frac{\mathrm{cond}(A,b)}{\mathrm{cond}(A)}\right)w_2\le \frac{\|x^*-x\|}{\|x^*\|}\le \mathrm{cond}(A)\left(1+\frac{\mathrm{cond}(A,b)}{\mathrm{cond}(A)}\right)w_2, 
\end{equation}
These bounds are more uniform in the sense that the direction of the right hand side has a smaller impact on the bounds. 
The relation between $w_1$ and $w_2$ is \citep{banoczi1998lack}
\begin{equation}\label{eq:w1_vs_w2}
    w_1 = \left(1 + \frac{\|A\|\|x\|}{\|b\|}\right)w_2.
\end{equation}
Note that $w_1$ is much larger than $w_2$ when $\|A\|\cdot\|x\| \gg \|b\|$ holds.
As a consequence, an approximation $x$ may yield $w_2$ that is smaller than $\tau$, thus satisfying criterion 2, but a large $w_1$, thus not satisfy criterion 1 for the same $\tau$, requiring more Krylov iterations to satisfy the criterion. On the other hand, $w_1$ has more chance to yield a solution of higher accuracy when $\mathrm{cond}(A,b)$ is small.
\subsection{Model order reduction}\label{s:mor}
Every optimization step requires the solution of two linear systems: one for the forward and one for the adjoint equation. We use model order reduction (MOR) to reduce the computational cost of these simulations, by constructing a surrogate model that provides
approximate solutions for the high-fidelity model. These solutions can be exploited as a solution
for the current linear system in order to avoid the CG solve completely
\citep{bennerMORsurvey}. 
In this paper, the surrogate model is built by projecting the high-dimensional discretized governing equations onto a low-dimensional reduced space. \\
Assume for simplicity that the objective is a given function of the discretized forward state variable $\mathbf{x}$. After discretization, equation \eqref{eq:energyEqn} results in
\begin{equation*}
    A(\rho^{\left(j\right)}) \mathbf{x}  =  b, \\
\end{equation*}
where $\rho^{\left(j\right)}$ is the design variable corresponding to the $j_{\text{th}}$ optimization iteration. Similarly, the discretized adjoint equation reads
\begin{equation*}
    A(\rho^{\left(j\right)}) \mathbf{y}  =  b_a, 
\end{equation*}
where $\mathbf{y}$ is the discretized adjoint state variable. The objective and the $i_{\text{th}}$ component of the gradient can be computed using
\begin{equation}\label{eq:discreteGrad}
    \left\{\begin{aligned}
        \mathcal{J}(\mathbf{x}(\rho^{\left(j\right)})) & =  f\left(\mathbf{x}\right), \\
        \frac{\partial \mathcal{J}}{\partial \rho_i} &= \mathbf{y}^T \frac{\partial
A}{\partial \rho_i} \mathbf{x}. 
    \end{aligned}\right.
\end{equation}
We build the reduced models by applying a Galerkin projection of the discretized forward and adjoint systems. Let $V$ be a $n\times r$ matrix with orthogonal columns spanning $r$ state vectors of the forward equations, i.e.,
\[
\mbox{Range}(V) = \mbox{span}\{ \mathbf{x}^{(f_1)}, \ldots,\mathbf{x}^{(f_r)}\},
\]
where $\mathbf{x}^{(j)} = \mathbf{x}{(\rho^{(j)})}$ the state vector of the forward equation corresponding to design variable $\rho$ at iteration $j$ of the optimization.
The forward reduced problem is
\[
    \begin{aligned}
            \widehat{A}(\rho) \widehat{\mathbf{x}} & =  \widehat{b}, \\
    \end{aligned}
\]
with $\widehat{A}(\rho) = V^T A(\rho) V$ and $\widehat{b} = V^T b$.\\
Let the columns of $W\in\mathbb{R}^{n\times s}$ span $s$ state vectors of the adjoint equation, i.e.,
\begin{equation*}
    \mbox{Range}(W) = \mbox{span}\{ \mathbf{y}^{(a_1)}, \ldots,\mathbf{y}^{(a_s)}\},
\end{equation*}
where $\mathbf{y}^{(j)} = \mathbf{y}{(\rho^{(j)})}$.
The adjoint reduced problem is
\[
    \widehat{A}(\rho) \widehat{\mathbf{y}} =  \widehat{b}_a, \\
\]
with $\widehat{A}(\rho) = W^T A(\rho) W$ and $\widehat{b}_a = W^T b$.\\
The solutions of both reduced systems i.e., $\widehat{\mathbf{x}}$ and $\widehat{\mathbf{y}}$, are projected back to the full order model space
\begin{align*}
    \widetilde{\mathbf{x}} &= V\widehat{\mathbf{x}}, \\
    \widetilde{\mathbf{y}} &= W\widehat{\mathbf{y}},
\end{align*}
where $\widetilde{\mathbf{x}}$ and $\widetilde{\mathbf{y}}$ are the MOR approximate solutions for the forward and the adjoint equations, respectively. The objective and the gradient can be approximated by using the equations
\begin{equation}\label{eq:redDiscreteGrad}
    \left\{\begin{aligned}
        \widetilde{\mathcal{J}}(\rho^{\left(j\right)}) & =  f\left(\widetilde{\mathbf{x}}\right), \\
        \frac{\partial \widetilde{\mathcal{J}}}{\partial \rho_i} &= \widetilde{\mathbf{y}}^T \frac{\partial
A}{\partial \rho_i}\widetilde{\mathbf{x}}. 
    \end{aligned}\right.
\end{equation}
It is important to note that when $r=s$ and the interpolation points of the forward reduced model coincide with the interpolation point of the adjoint reduced model i.e., $\rho^{(f_i)}=\rho^{(a_i)}$ for every $i$, then the objective and the gradient are interpolated exactly at the interpolation points:
\begin{eqnarray*}
    \left\{\begin{aligned}
            \mathcal{J}(\rho^{(j)}) & =  \widetilde{\mathcal{J}}(\rho^{(j)}),\quad j=1,\ldots,r \\
            \nabla_\rho \mathcal{J}(\rho^{(j)}) & =  \nabla_\rho
            {\widetilde{\mathcal{J}}}(\rho^{(j)}),
    \end{aligned}
    \right.
\end{eqnarray*}
 When the reduced models
are updated at different design values $\rho$, the gradient is not
interpolated exactly at the interpolation points but it is well approximated whenever
the state vectors are well approximated.
In particular, the error on the $i_{th}$ component of the gradient computed by the reduced model at interpolation points $\rho^{(j)}$ is given by
\begin{equation*}
    E_{\nabla\mathcal{J}_i} = {\Delta\mathbf{ y}}^T \frac{\partial A}{\partial \rho_i} {\mathbf{x}} + {\mathbf{ y}}^T \frac{\partial A}{\partial \rho_i} {\Delta\mathbf{  x}} + {\Delta \mathbf{y}}^T \frac{\partial A}{\partial \rho_i} {\Delta \mathbf{x}},
\end{equation*}
where $\Delta\mathbf{ x}= \widetilde{\mathbf{x}} - \mathbf{x}$ and $\Delta\mathbf{
y}= \widetilde{\mathbf{y}} - \mathbf{y}$ are the state errors of the forward and the
adjoint MOR, respectively. These errors can be bound using the corresponding MOR residuals using equation \eqref{eq:bound_stop_1} or \eqref{eq:bound_stop_2} depending on the criterion used. Note that the forward state error is zero at every forward interpolation point, i.e., $\Delta\mathbf{x}=0 \text{ at } \rho=\rho^{(f_1)},\ldots, \rho^{(f_r)}$ and similarly for the adjoint state error. \medskip \\
An advantage of keeping two different reduced models is that $V$ and $W$ do not need to have the same number of columns and can be updated independently. This leads to a more effective use of the reduced models. For instance, if the forward field can be better approximated than the adjoint field, $V$ can be updated less frequently, reducing computational cost.
The remaining question is how to compute the reduced basis $V$ and $W$ out of $r$ (or
$s$) high-fidelity simulations. We choose to simply orthogonalize these vectors using
a Gram-Schmidt procedure, following
\citep{GoguMorTopopt,choi2019acceleratingdesignoptimizationusing}. To keep the size of
the bases small, we  only keep the $r$ latest interpolation points, 
throwing away the older points. We progressively update the bases as new vectors
become available, as it is described in \ref{alg:gso}.
\begin{algorithm}[tp]
    \begin{algorithmic}[1]
        \STATE Given vector $x \in \real^{m}$ 
        \STATE Given base $V \in \real^{m \times k}$
        \STATE Given base $R \in \real^{k \times k}$
        \IF{$k \le r$}
          \STATE $[v,r] \gets \mathrm{GS}\left(V,x\right)$ \COMMENT{Orthogonalize x
          against V using Gram-Schmidt (GS)}
          \STATE $V \gets \left[V, v\right]$ \COMMENT{Append new vector to $V$}
          \STATE $R \gets \left[R, r\right]$ \COMMENT{Append new vector to $R$}
    \ELSE
            \STATE $V \gets \left[V, x\right]$ \COMMENT{Append new vector to $V$}
            \STATE $[Q, R] \gets \mathrm{QR}\left(V\right)$ \COMMENT{Orthogonalize
            $V$ using GS}
            \STATE $[Z, S] \gets \mathrm{QR}\left(R(:,2:r+1)\right)$ \COMMENT{Orthogonalize the last $r$ columns of R}
            \STATE $V \gets QZ$
    \ENDIF	
    \RETURN { Q and R }
        \end{algorithmic}
    \caption{Procedure to add a new vector to the base vector.}
    \label{alg:gso}
\end{algorithm}
\subsection*{Model order reduction workflow}
\begin{figure}[tp]
	\centering
	\input{mor_workflow_diagram}
        \caption{Incorporating MOR techniques into the topology optimization workflow.}\label{fig:morFlowchart}
\end{figure}
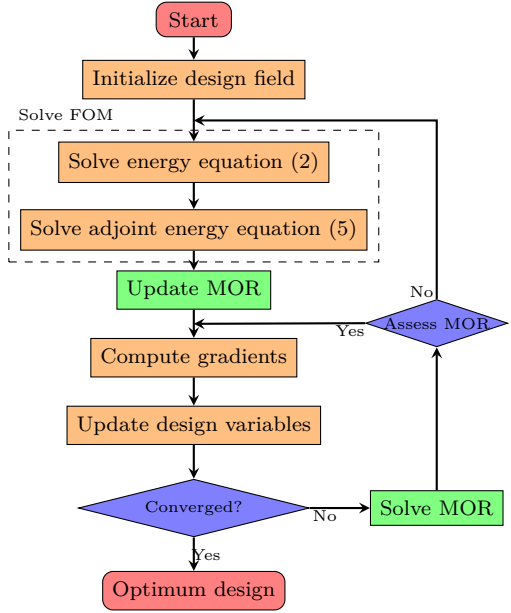

Figure~\ref{fig:morFlowchart} summarizes the proposed workflow for integrating MOR
into the topology optimization process described on Figure~\ref{fig:normalFlow}. Similar to the standard approach,
the optimization begins with initializing the design field, followed by the FOM
analysis that consists in solving both forward and adjoint equations.
The resulting FOM solutions are then used to expand the reduced bases in the step
\textbf{update MOR}, using 
\ref{alg:gso}.
For every new design, approximations are first computed using the reduced models in
the \textbf{solve MOR} step. The accuracy of MOR is then assessed using the same
criterion \eqref{eq:stop2} used by the conjugate gradient solver, but with a reduced threshold
$\romtau$. 
This evaluation step is referred to as \textbf{assess MOR}. If
the MOR approximations meet the stop criterion with reduced threshold $\romtau$, they are accepted and gradients are computed directly from them. Otherwise, the workflow solves the FOM using the conjugate gradient method until the stop criterion is satisfied with threshold $\tau$. \\ 
This switch between MOR and FOM is driven solely by the residual. The FOM solver is only triggered when the MOR approximation fails to satisfy the residual tolerance $\romtau$. This adaptive strategy allows expensive FOM solves to be performed only when the surrogate model is not accurate enough anymore, improving computational efficiency.\\
If either the state or the adjoint model fails to satisfy the tolerance criterion, only
the corresponding equation is solved up to convergence using the conjugate gradient
solver, while we still rely on the MOR approximation for the other equation. 
\section{Numerical experiments}\label{s:experiments}
In this section, we investigate the efficiency of the proposed workflow for 2D and 3D topology optimization benchmark problems. First, in section \ref{ss:stopping} we show the importance of choosing the stopping criterion \eqref{eq:stop2} instead of \eqref{eq:stop1} on the
performance of MOR for a 2D problem. Second, in section \ref{ss:MORaccurate} we demonstrate the efficiency for a 3D problem where we combine MOR with both, accurate and partially inaccurate, FOM solvers. \\
For the optimization, we rely on a C++ implementation of the method of moving assymptotes (MMA) method developed in
\citep{aage2015topology,svanberg1987method}. 
Furthermore, we use a constant penalization exponent equal to $3$ in the SIMP material interpolation \eqref{eq:SIMP}.\\ 
For the simulation of the forward and the adjoint equations, we rely on OpenFOAM version 2306 where the Finite Volume Method (FVM) is applied on a mesh with rectangular elements to discretize the governing forward and adjoint equations \eqref{eq:energyEqn} and \eqref{eq:adjEnergyEqn}, as well as the PDE filter equation \eqref{eq:pde_filter}. The unknown variables i.e., density field and state variables, are piecewise constant functions characterized by their values at the center of every mesh cell. The resulting linear systems are solved using a preconditioned conjugate gradient relying on the OpenFOAM geometric agglomerated algebraic multigrid (GAMG) implementation for the preconditioner with the settings shown in Table \ref{tab:GAMG_config}. Note that the preconditioner uses a geometric agglomeration that depends on the mesh. Since the mesh does not change along the optimization, the agglomeration is computed once and saved in memory to avoid unnecessary computations.
\begin{table}[ht]
    \centering
    \begin{tabular}{|c|c|}
        \hline
        \textbf{Name}     & \textbf{Value} \\ \hline
         smoother             & symGaussSeidel \\ \hline
         nVcycles             & 1 \\ \hline
         maxPreSweeps         & 1 \\ \hline
         maxPostSweeps        & 1 \\ \hline
         nPreSweeps           & 1 \\ \hline
         nPostSweeps          & 1 \\ \hline
         directSolveCoarsest  & yes \\ \hline
         cacheAgglomeration   & yes \\ \hline
    \end{tabular}
    \caption{GAMG OpenFOAM configuration settings.}
    \label{tab:GAMG_config}
\end{table}
Finally, to solve the reduced model, we simply use a direct solver based on LU factorization implemented in OpenFOAM.
\subsection{Importance of the stopping criteria}\label{ss:stopping}
\subsubsection{Problem description} \label{subsub:problem_descrip_2D}
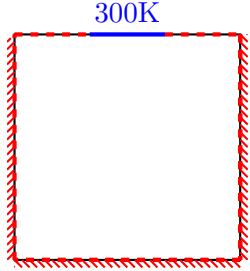
\begin{figure}[tp]
\begin{center}
    \input{figure_geometry}
	\caption{Geometry and boundary conditions of the problem setup.}\label{fig:caseGeometry}
\end{center}
\end{figure}

\begin{figure}[tp]
	\centering
	\begin{subfigure}[b]{0.45\textwidth}
		\centering
		\includegraphics[width=\textwidth]{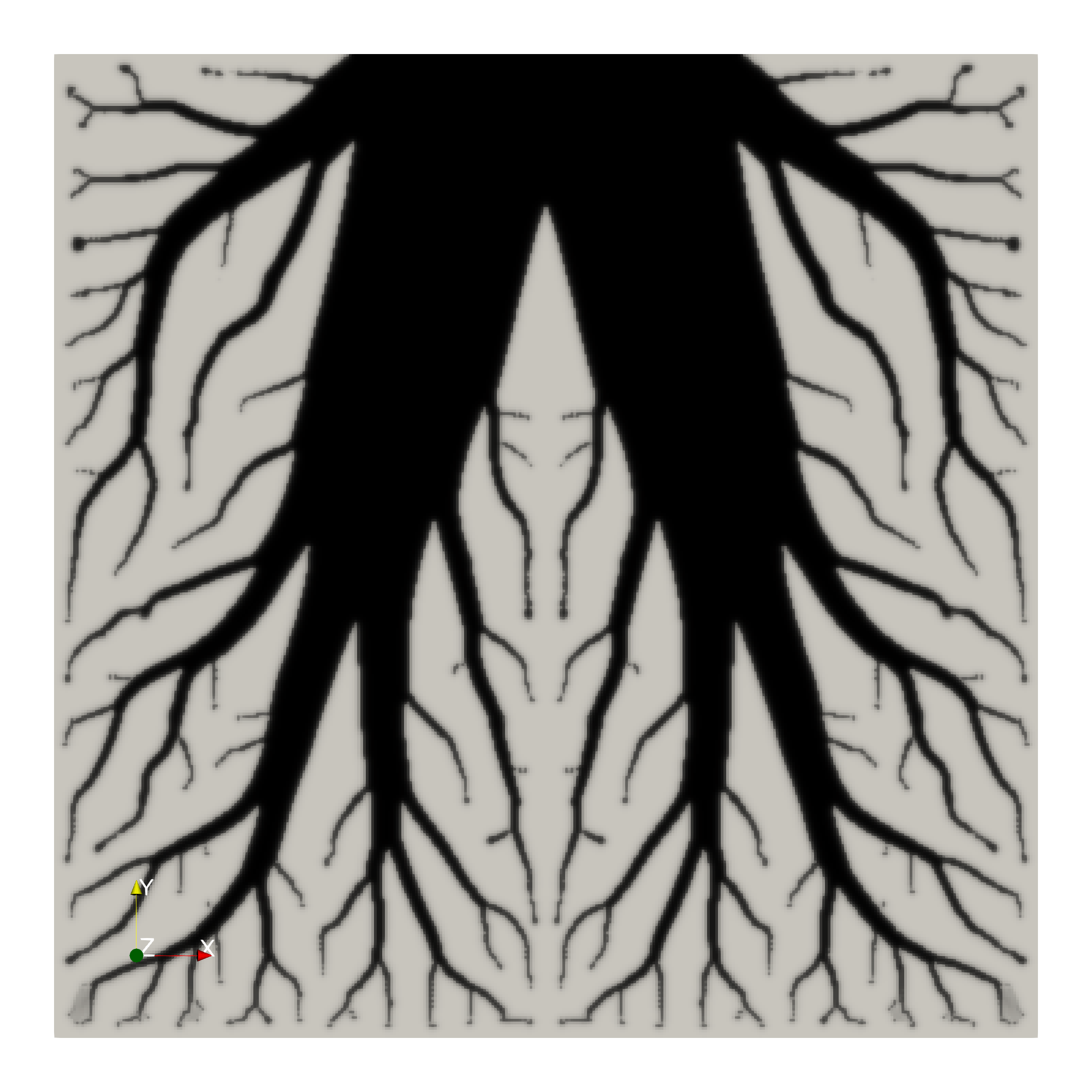}
		\caption{Stop criterion defined by \eqref{eq:stop1} with $\tau=10^{-13}$ and $\romtau =5\cdot 10^{-3}$}
            \label{subfig:eps1_5e-3}
	\end{subfigure}
	\hfill
	\begin{subfigure}[b]{0.45\textwidth}
		\centering
		\includegraphics[width=\textwidth]{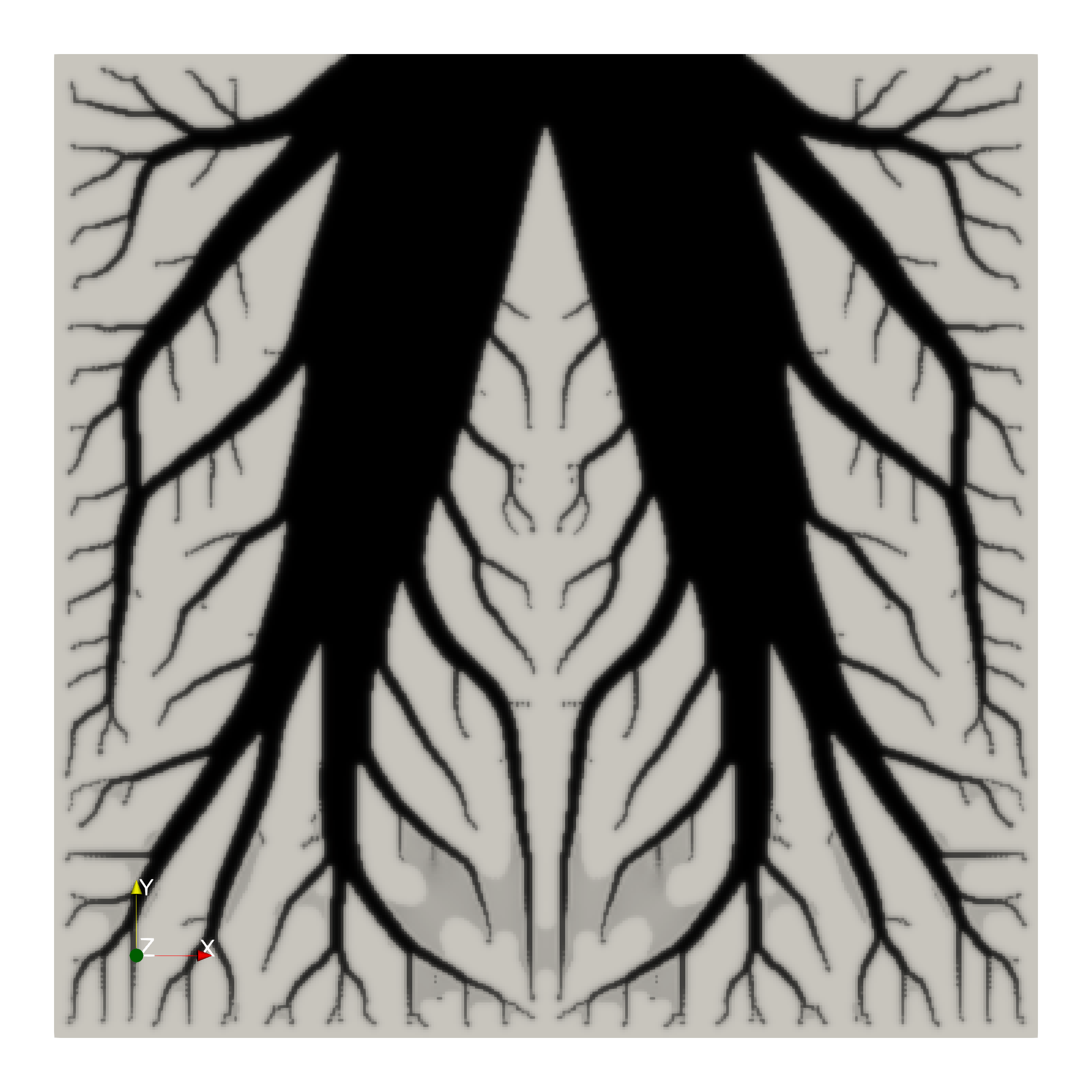}
		\caption{Stop criterion defined by \eqref{eq:stop2} with $\tau=10^{-13}$ and $\romtau=10^{-6}$}
        \label{subfig:eps2_1e-6}
	\end{subfigure}
	\hfill
    \caption{Density field $\rho$ at iteration $500$ for the 2D problem of subsection \ref{ss:stopping}. Black is the high conductive material and grey is the low conductive material.}
    \label{fig:2D-designs}
\end{figure}

\begin{figure}[tp]\label{fig:comparisonResiduals}
	\centering
	\begin{subfigure}[b]{0.45\textwidth}
		\centering
		\includegraphics[width=\textwidth]{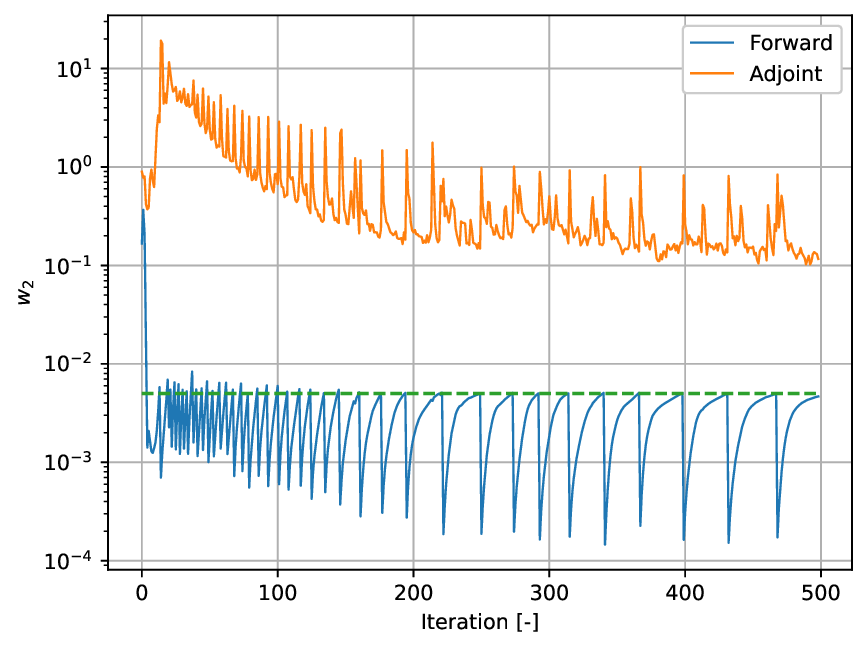}
		\caption{$w_1$ values as defined in \eqref{eq:stop1}, $\tau=10^{-13}$ and $\romtau =5\cdot 10^{-3}$}
        \label{subfig:res1_5e-3}
	\end{subfigure}
	\hfill    	
    \begin{subfigure}[b]{0.45\textwidth}
		\centering
        \includegraphics[width=\textwidth]{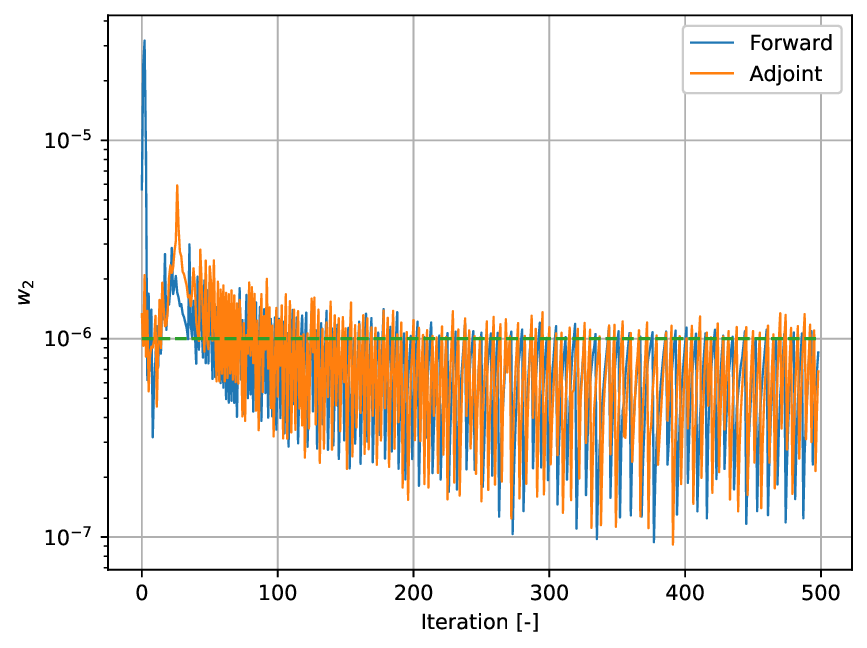}		\caption{Stopping criterion \eqref{eq:stop2}, $\tau=10^{-13}$ and $\romtau=10^{-6}$}
        \label{subfig:res2_1e-6}
	\end{subfigure}
    \hfill
    \caption{Computed MOR residuals as a function of the design iterations for the
    forward and the adjoint equations for the 2D problem of subsection \ref{subsub:problem_descrip_2D}. Dashed line is the corresponding $\romtau$ value.}
\end{figure}

We investigate the influence of the stopping criterion on the overall performance of the
reduced model.  To this end, we consider the setup of \citep{gersborg2006topologyHeat} whereby the domain $\Omega$, shown in Figure~\ref{fig:caseGeometry}, consists of a square domain of size $12\times12$ where a fixed temperature is applied to the middle part of the top boundary temperature
and the remaining part of the boundary acts as an insulated wall. Additionally, a constant heat flux of $q_n=1000$ is applied in the out-of-plane direction. The domain is discretized using $360$ elements on each side resulting with a total of $129 600$ mesh cells. \\
In order to make the problem not self-adjoint, we choose the average square shifted temperature 
over the domain as the objective function
 \begin{equation*}
      \mathcal{J}(\enthalpy) = \frac{1}{|\Omega|}\int_\Omega \left(\enthalpy - 300\right)^2 d\Omega,
 \end{equation*}
with a volume fraction constraint equal to $v_f=0.4$ of the domain volume.  The optimization is initialized with a uniform design satisfying the volume constraint i.e., $\rho^{(0)}=0.4 \text{ in } \Omega$. Moreover, we set the MMA move limit to $0.1$ and we use \eqref{eq:relation_filter} with $\lambda_h$ twice the mesh size to set the length parameter of the PDE regularizing filter \eqref{eq:pde_filter} to  $\lambda=0.028$.
The residual threshold used to stop the linear solver is set to $\tau = 10^{-13}$ while the reduced threshold $\romtau$ is selected based on the chosen stop criterion reported in table \ref{tab:solver-config}. The maximum size of the MOR basis is set to $r=10$.
\begin{table}[ht]
	\centering
	\caption{Solver configuration}
	\label{tab:solver-config}
	\begin{tabular}{|l|c|c|c|}
		\hline
		\textbf{Criterion} & $\tau$ (FOM) & $\romtau$ (MOR)&r\\
		\hline
		\hline
		  $\|r^{(j)}\| \leq \tau \|b\|$  & \(10^{-13}\) &$5\cdot 10^{-3}$&10 \\
		\hline
		$\|r^{(j)}\| \leq \tau (\|b\| + \|A\| \|x^{(j)}\|)$& \(10^{-13}\)  & $10^{-6}$ &10 \\
		\hline
	\end{tabular}
\end{table}
Figure \ref{subfig:eps1_5e-3} and \ref{subfig:eps2_1e-6} show the final designs resulting from the optimization process using the two different stopping criteria \eqref{eq:stop1} and \eqref{eq:stop2} respectively. Both designs exhibit the classical tree-like structures. Figures~\ref{subfig:res1_5e-3} displays the resulting error measure $w_1$ values, as defined in \eqref{eq:stop1}, when computed using MOR solutions. Similarly, \ref{subfig:res2_1e-6} shows the resulting $w_2$ values as defined in \eqref{eq:stop2}. 
When using stopping criterion \eqref{eq:stop1}, a discrepancy arises between the performance of the forward and adjoint MORs. Specifically, the adjoint MOR solutions do not satisfy the stopping criterion since the computed $w_1$ values are above the threshold (constant green line) leading to a poor use of the adjoint reduced model. This is in contrast to the forward MOR, which is utilized extensively throughout the design process.
In constrast, when using criterion \eqref{eq:stop2} instead, we observe that the reduced models exhibit a similar performance because the $w_2$ error measure are of the same order of magnitude and below the threshold. This results in a more effective use of the reduced model because the computed approximations are used more often during optimization. \\ 
As discussed in section \ref{ss:stop}, the relation between $w_1$ and $w_2$, given by equation \eqref{eq:w1_vs_w2}, is determined by the ratio $\nicefrac{\|A\| \|x\|}{\|b\| }$.
Table~\ref{tab:norms} displays the norm of the matrix $A$ multiplied by the norm of the state vector $x$, of the source $b$ and the corresponding ratio $\nicefrac{\|A\| \|x\|}{\|b\| }$  computed at the final iteration of the design process. From the table, we can see that for both, the forward and adjoint equations, the inequality \( \|A\| \|x\| \gg \|b\| \) holds.
However, this inequality is much more pronounced for the adjoint equation by a factor of $10^6$ compared to the $10^2$ for the forward equation.
This explains the discrepancy observed in the figure between the forward and the adjoint reduced model when the stopping criterion \eqref{eq:stop1} is used, where the behavior is dependent on the right hand side. In principle, we could use a different threshold value $\romtau$ for the forward and the adjoint equations, but this would introduce an extra layer of tuning that is not desirable.
The stopping criterion \eqref{eq:stop2} seems consistently less dependent on the right hand side.  In the remainder of the paper, we adopt the stopping criterion defined by equation \eqref{eq:stop2}. 
\begin{table}[tp]
    \centering
    \begin{tabular}{|l|c|c|c|}
        \hline
        \textbf{}                 & $\norm{A}\norm{x} $    & $\norm{b}$ & $\frac{\norm{A}\norm{x}}{\norm{b}}$      \\ \hline
        Forward equation          & $6.39\cdot 10^{8}$     &  $1.97\cdot 10^{6}$ & $3.25\cdot 10^{2}$     \\ \hline
          Adjoint equation          & $5.07\cdot 10^6$    & $4.56$ &  $1.11\cdot 10^6$             \\ \hline
    \end{tabular}
    \caption{Norms of $Ax$, the right-hand side $b$ and their ratio \eqref{eq:w1_vs_w2} for the forward and the adjoint equation at the final iteration of the optimization for the 2D problem of subsection \ref{subsub:problem_descrip_2D}.}\label{tab:norms}
\end{table}
\subsection{Comparison of MOR against FOM solvers}\label{ss:MORaccurate}
This section aims to assess the computational savings of MOR in the proposed workflow for a 3D
thermal topology optimization problem by comparison with alternative solution strategies for the state and adjoint equations based on the multigrid preconditioned CG method. 
First, we evaluate the speedup of MOR against the standard strategy which solves every linear system up to full accuracy.
However, we find that despite achieving a significant speedup, the method remains overly less performant than the recently proposed one-shot method of \citep{amir2024one} that consists in solving every linear system using only one iteration of the multigrid preconditioned conjugate gradient (MGCG) method. Therefore, we propose an alternative strategy where we combine MOR with the one-shot method by constructing the reduced basis from one-shot solutions. We find that an even higher speed up can be achieved by this proposed solver strategy. Note that because the  focus lies on the forward and adjoint equations, the computational time associated with the PDE filtering step is excluded from the reported timings. Nevertheless, the filter equation could be accelerated efficiently using the same MOR methodology.
\subsubsection{3D problem setup}
\begin{figure}[tp]
\begin{center}
    \input{3D_geometry}
	\caption{Geometry and boundary conditions of the 3D problem setup.}\label{fig:3DcaseGeometry}
\end{center}
\end{figure}
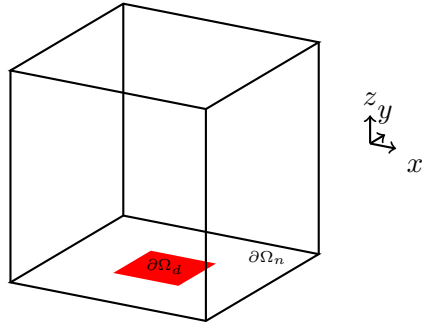

 We consider the design problem formulated in \citep{FEPPON2020109574} and illustrated on figure  \ref{fig:3DcaseGeometry}. An uniform heat source is applied to a cubical domain where a Dirichlet boundary
 condition is imposed on the red square $\partial \Omega_d$ with a temperature of
 $T_d=273$.
 On the remaining boundaries $\partial \Omega_n$, a Neumann boundary condition
 is modelling insulating walls. This domain is discretized into a $200\times200\times200$ structured mesh in $8$ million mesh quadrilateral elements.
 Additionally, the whole domain is heated with a spatially constant source $Q=10^4$. The goal is to design a heat
 sink by optimally distributing two materials of conductivities $\kappa_s=1$ and $\kappa_d=100$ throughout the domain $\Omega$ subject to a volume constraint of $v_f=0.05$. 
Again to avoid self-adjointness, we choose  the average square shifted temperature over the domain as the objective function:
 \begin{equation*}
      \mathcal{J}(\enthalpy) = \frac{1}{|
      \Omega|}\int_\Omega \left(\enthalpy-273\right)^2 d\Omega.
 \end{equation*}
 Similarly to previous section, the optimization is initialized with an uniform design with $\rho=0.05$. The MMA \textbf{move\_limit} variable is set to $0.1$ and we use \eqref{eq:relation_filter} with $\lambda_h$ twice the mesh size to set the length parameter of the PDE regularizing filter \eqref{eq:pde_filter} to  $\lambda=0.0028$.
\subsubsection{MOR against CG until convergence}
\begin{figure}[H]
	\centering
    \begin{subfigure}[b]{0.45\textwidth}
		\centering
        \includegraphics[width=\textwidth]{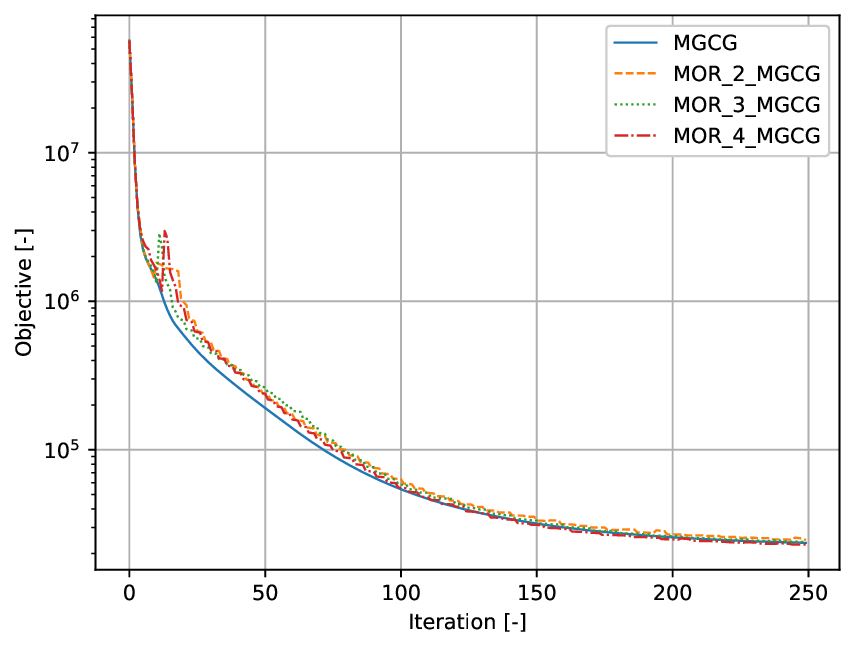}
        \caption{Objective value}
	\end{subfigure}
    \begin{subfigure}[b]{0.45\textwidth}
		\centering
        \includegraphics[width=\textwidth]{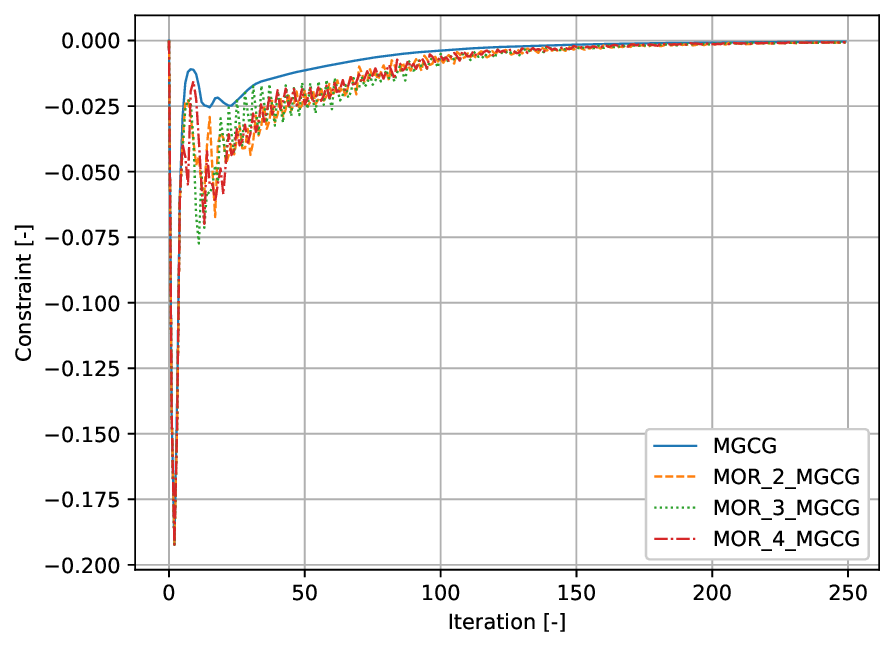}
        \caption{Constraint value}
	\end{subfigure}
    \caption{Convergence history of the objective and constraint values for the 3D thermal optimization problem using MOR in combination of different accurate FOM solvers.}
	\label{fig:MOR_full_iter_obj_cons}
\end{figure}

\begin{figure}[H]
	\centering
	\begin{subfigure}[b]{0.24\textwidth}
		\centering
		\includegraphics[width=\textwidth]{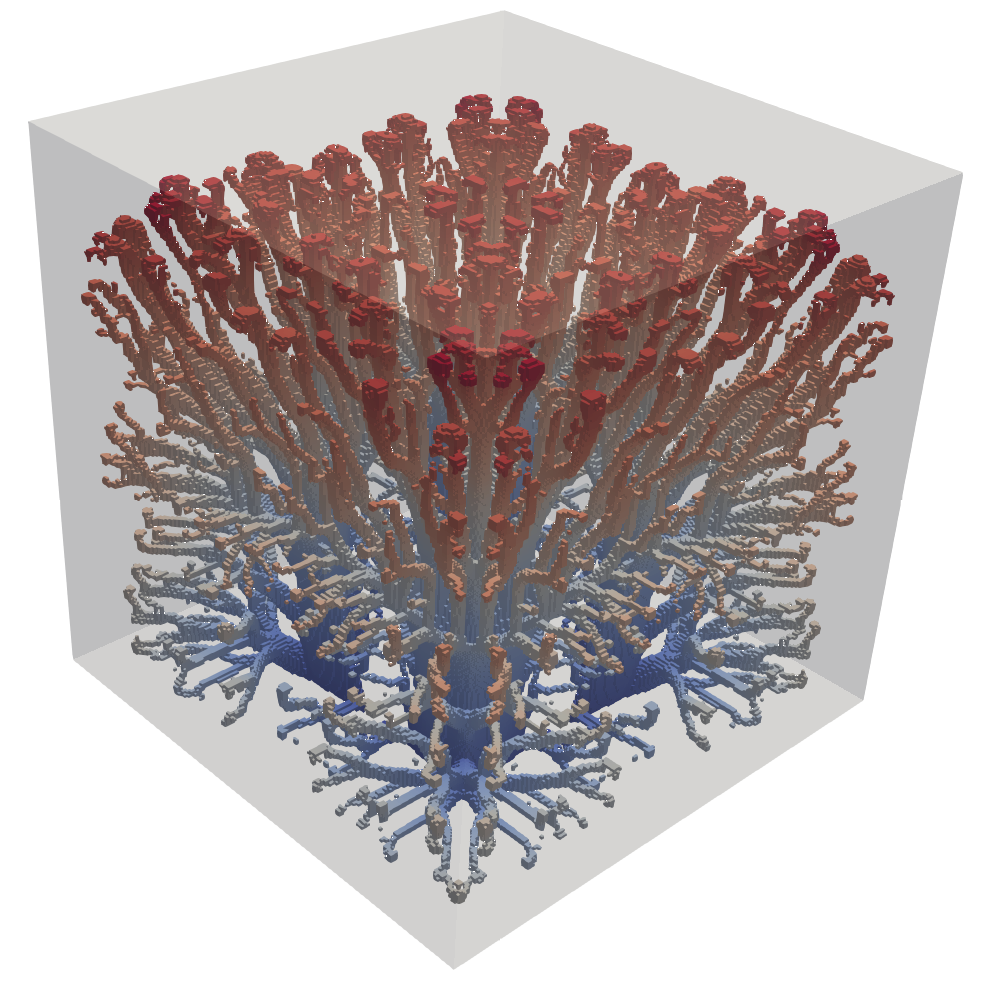}
		\caption{MGCG}\label{subfig:acc:MGCG_full_iter}
	\end{subfigure}	
	\begin{subfigure}[b]{0.24\textwidth}
		\centering
		\includegraphics[width=\textwidth]{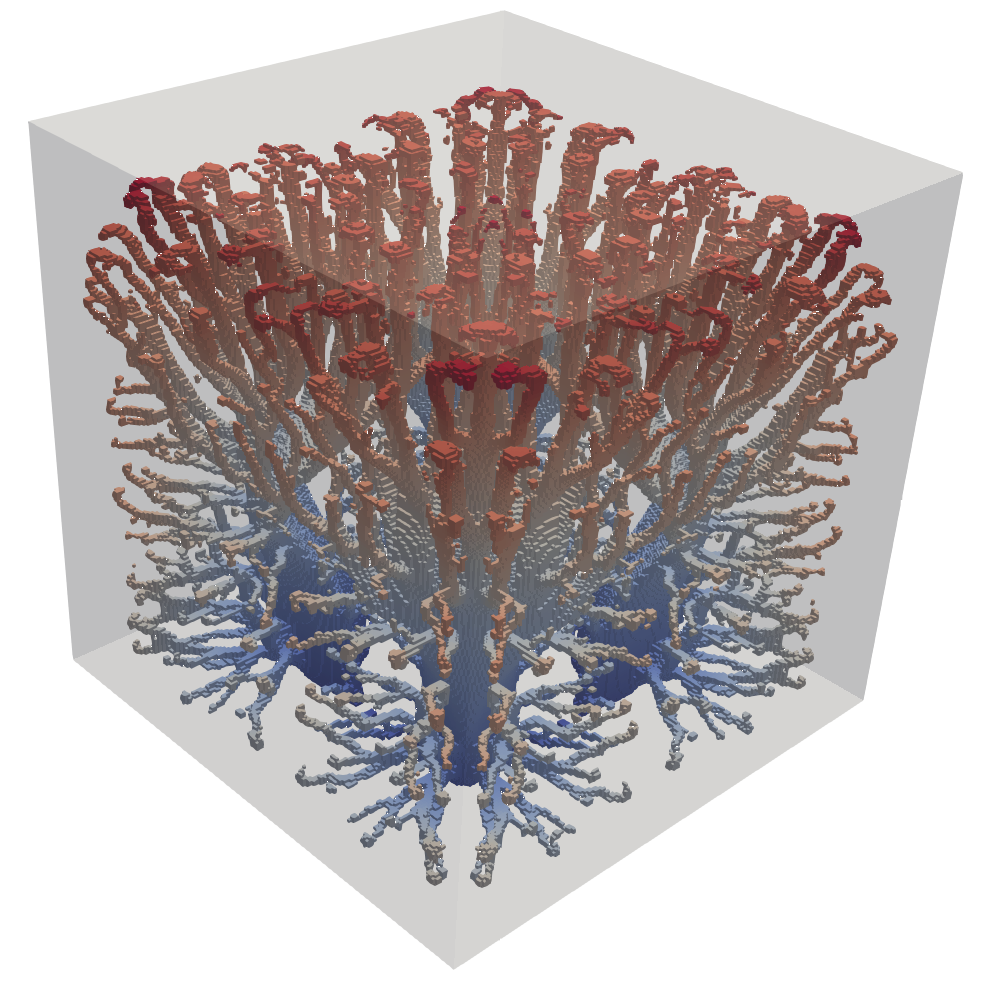}
		\caption{MOR\_2\_MGCG}\label{subfig:acc:MOR_2_MGCG_2e-6_full_iter}
	\end{subfigure}
	\begin{subfigure}[b]{0.24\textwidth}
		\centering
		\includegraphics[width=\textwidth]{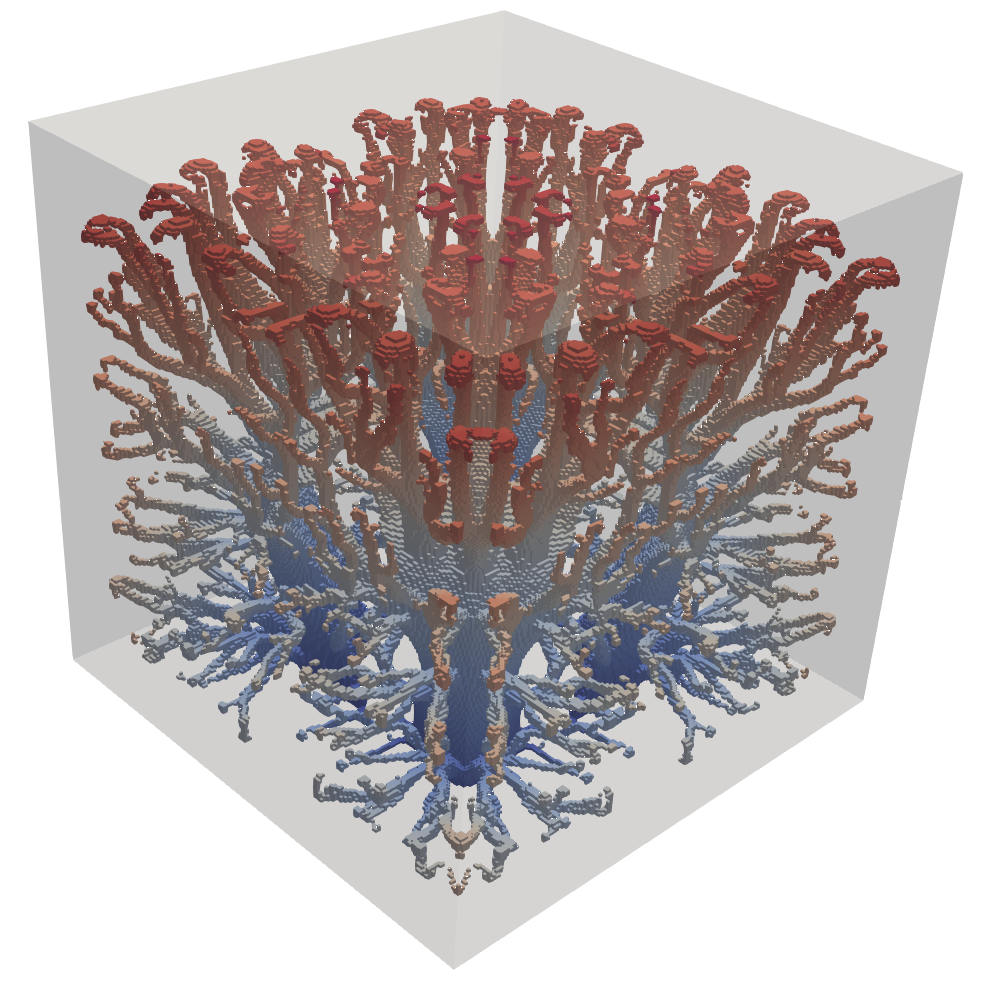}
	   \caption{MOR\_3\_MGCG}\label{subfig:acc:MOR_3_MGCG_2e-6_full_iter}
	\end{subfigure}
	\begin{subfigure}[b]{0.24\textwidth}
		\centering
		\includegraphics[width=\textwidth]{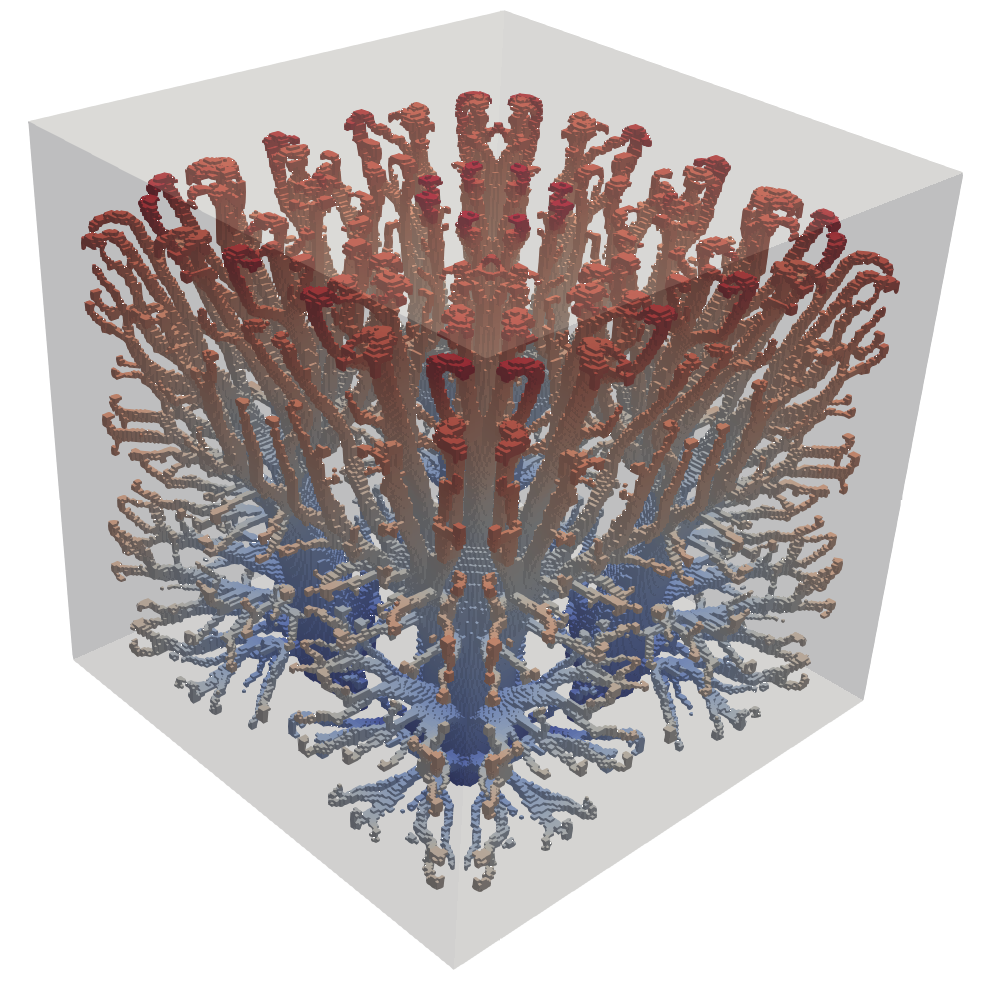}
		\caption{MOR\_4\_MGCG}\label{subfig:acc:MOR_4_MGCG_2e-6_full_iter}
	\end{subfigure}
	\caption{Optimal designs of the 3D thermal optimization using MOR in combination of different accurate FOM solvers. (a) FOM until convergence, (b) MOR with $r=2$, (c) MOR with $r=3$ and (d) MOR with $r=4$ in combination with FOM until convergence. A clip value of $0.5$ is used for the visualization}
	\label{fig:MOR_full_iter_eps}
\end{figure}

\begin{figure}[H]
	\centering
    \includegraphics[width=0.4\textwidth]{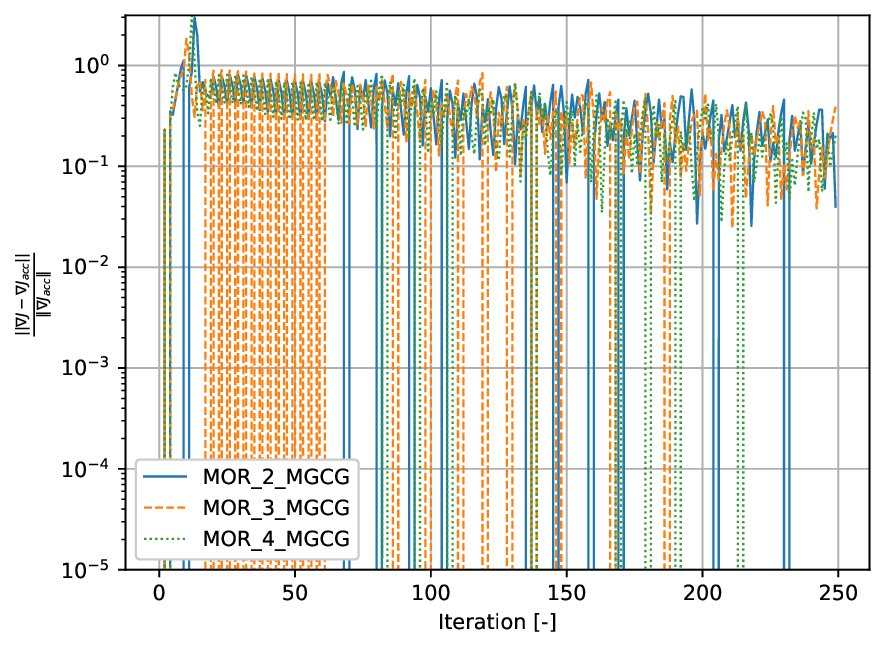}
    \caption{Relative error of the MOR approximation for gradient with respect to the gradient with respect to the accurate gradient i.e., computed using fully converged forward and adjoint solutions.}    
    \label{fig:3D_1_ang}
\end{figure}

\begin{figure}[H]
	\centering
    \includegraphics[width=0.4\textwidth]{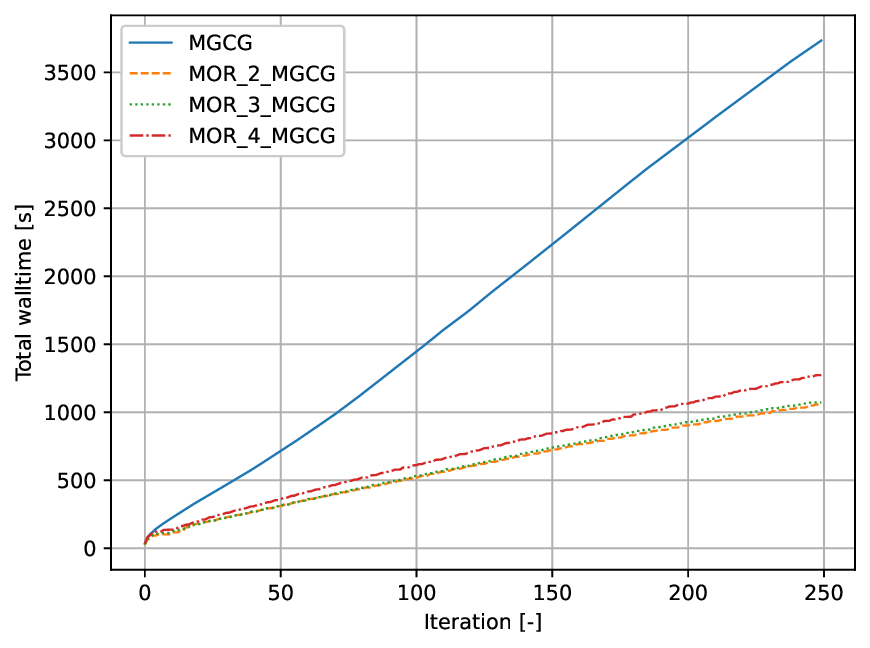}
    \caption{Total walltime as a function of the iteration number for the workflow with only accurate FOM solver and the workflow with MOR incorporated for $250$ iterations.}    
    \label{fig:acc:timings}
\end{figure}

In this section, we investigate the performance of MOR when we solve the state and adjoint FOM equations \eqref{eq:energyEqn} and \eqref{eq:adjEnergyEqn} until convergence. To this end, the reduced order model is constructed from solutions from MGCG until convergence. We test out methodology with varying dimensions of the reduced basis .The different combinations, the number of basis vectors and the corresponding stopping criterion thresholds used in this experiment are reported in table \ref{tab:solver-config-3D-1}. \\
\begin{table}[ht]
	\centering
	\begin{tabular}{|l|c|c|c|}
		\hline
		\textbf{Name} &$r$ &$\fomtau$& $\romtau$ \\
		\hline
        MGCG &NA&$1 \cdot 10^{-13}$& NA \\
		\hline
		MOR\_2\_MGCG & 2&$1 \cdot 10^{-13}$& $5\cdot 10^{-6}$ \\
		\hline        		
        MOR\_3\_MGCG & 3&$1 \cdot 10^{-13}$& $5\cdot 10^{-6}$ \\ \hline
        MOR\_4\_MGCG & 4&$1 \cdot 10^{-13}$& $5\cdot 10^{-6}$ \\ \hline
	\end{tabular}
    \caption{Solver configuration for the 3D thermal optimization using MOR in combination of different accurate FOM solvers. $r$ is the maximum dimension of the reduced bases, $\fomtau$ and $\romtau$ are the FOM and MOR residual threshold respectively. }
	\label{tab:solver-config-3D-1}
\end{table}
Figure \ref{fig:MOR_full_iter_obj_cons} displays the objective and the constraint values along the optimization path for the different solver strategies. We can see that the 4 solver strategies have similar objective and constraint values during the optimization process decreasing and flattening out around iteration $200$. Interestingly, the three MOR strategies have a slower convergence in the first $60$ iterations but reach almost identical performances to that of the full order model solver, at convergence. Similarly, the constraint value is better enforced when using the high fidelity MGCG solver in the first $100$ iterations. However, there are no visible differences after that point. Note that around iteration 12, there is a visible spike in the objective function for all MOR solver strategies. This occurs when the forward MOR approximation fails to meet the required residual threshold after several consecutive MOR steps, triggering a fallback to the accurate MGCG solver, that corrects the accumulated error in the temperature field, resulting in an adjustment in the objective value. The smaller magnitude of subsequent corrections can be attributed to the smaller design changes as the optimization progresses. These smaller changes lead to MOR approximations that remain closer to the exact converged fields, thereby minimizing the need for large adjustments.
\\
Figure \ref{fig:MOR_full_iter_eps} shows the  final designs generated by the different solution strategies, with and without MOR integration. Figure \ref{subfig:acc:MGCG_full_iter} depicts the accurate strategy with only the FOM solver until convergence while Figures \ref{subfig:acc:MOR_2_MGCG_2e-6_full_iter}, \ref{subfig:acc:MOR_3_MGCG_2e-6_full_iter} and \ref{subfig:acc:MOR_4_MGCG_2e-6_full_iter} illustrate the resulting optimized designs when the MOR solver is used. The optimized material distribution in each scenario predominantly places the highly conductive material adjacent to the fixed temperature boundary at the bottom, which produces the highest heat flux. It is worth noting that we do not aim to strictly reproduce the exact design of the reference MGCG accurate solver. As shown in Table \ref{tab:MGMOR_full_iter}, the optimized designs using MOR achieve a performance within $4\%$ of the one obtained with accurately solving the state and adjoint equation are every iteration. 
To quantitatively evaluate the accuracy of the MOR approach, we compute the relative error of the gradient as
\begin{equation}
\frac{\|\nabla \mathcal{J}-\nabla \mathcal{J}_{\mathrm{acc}}\|}{\|\nabla \mathcal{J}_{\mathrm{acc}}\|},    
\end{equation}
where $\|\nabla \mathcal{J}_{acc}\|$ denotes the accurately computed gradient obtained by fully converged forward and adjoint solutions. 
Figure \ref{fig:3D_1_ang} plots this metric for every iteration, showing an oscillatory pattern.
This pattern arises from periodic updates of the reduced bases when approximations do not meet the specified accuracy threshold. The error  small at iterations where both bases, forward and adjoint, are updated simultaneously.
Note that the error in the gradient norm is large at the first optimization iterations.
However, these inaccuracies do not seem to cause instabilities that break the optimization and could potentially be mitigated by setting a smaller reduced residual threshold $\romtau$ to improve the accuracy of the MOR approximations.\medskip
Table \ref{tab:MGMOR_full_iter} summarizes the computational performance and the final objective function values obtained. The table reports the walltime, defined as the total time spent in solving the linear systems during the whole optimization. The MOR solution strategies achieve speedups of $3.50$ while still producing the optimized designs wit performance is falling within $4\%$ of the strategy based on the full order model. The computational speedup peaks at $r=2$ (3.50x) but it is still significant for $r=4$ (2.96x). Interestingly, for this problem, increasing the basis size has a small impact on the total number of reductions which suggests that the smaller speedup can be attributed to the extra overhead of maintaining the larger bases. In other words, the computational overhead of managing and projecting a larger reduced basis outweighs the marginal reduction in solver iterations achieved by the richer subspace. However the final optimized design is of slightly better quality regarding the performance and the satisfaction of the constraints. To provide context for these timings, the time spent on other operation of the optimization procedure (such as computing the design updates) was approximately $2065.82$s, indicating that in the reference case, the linear solver remains a significant portion of the total computational cost. 
Finally, Figure \ref{fig:acc:timings} shows the total walltime as a function of the iteration number where it is clear that the MOR solver strategies achieve faster results already at the very beginning of the optimization. Moreover, from the flattening of the MOR curves, it is clear that MOR is exploited even more towards the end of the optimization, where the optimized designs stagnates.
\begin{table}[ht]
  \centering
  \begin{tabular}{|l|c|c|c|c|c|}
    \hline
                    & Objective & Walltime  & Speedup   & Forward    & Adjoint    \\ 
                    & Value     & (s)       & (CG)       & Reductions & Reductions \\ \hline
                    MGCG & 23620.30 &  3722.67  &  1 & N/A & N/A \\ \hline
                    MOR\_2\_MGCG & 24616.06 &  1061.55  &  3.50& 189 & 176 \\ \hline
                    MOR\_3\_MGCG & 23859.72 &  1070.53  &  3.47& 188 & 175 \\ \hline
                    MOR\_4\_MGCG & 22820.57 &  1269.48  &  2.96& 189 & 176 \\ \hline
  \end{tabular}
  \caption{Numerical results of 3D thermal topology optimization using MOR in combination with MGCG until convergence.}\label{tab:MGMOR_full_iter}
\end{table}
\subsubsection{MOR against the one-shot method}
\begin{figure}[h]
	\centering
    \includegraphics[width=0.40\textwidth]{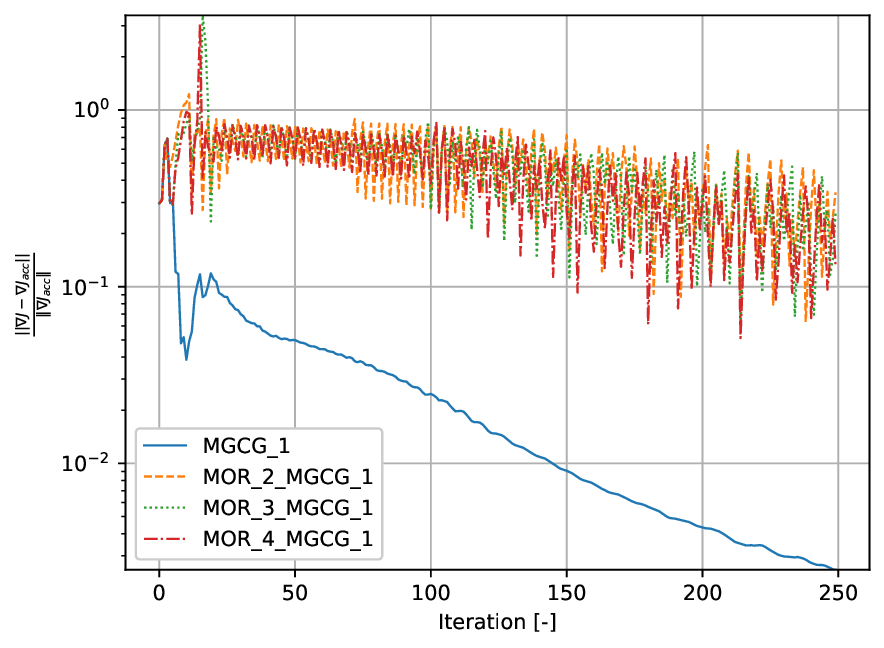}
    \caption{Relative error of the computed gradient with respect to the gradient computed with the accurate gradient, i.e., computed with forward and adjoint equations solved accurately until convergence. FOM solver is MGCG.}    \label{fig:3D_2_err}
\end{figure}

The one-shot from \citep{amir2024one} recently arised in topology optimization. It consists applying one iteration of the multigrid preconditioned conjugate gradient method to every linear systems from the forward and the adjoint equations \eqref{eq:energyEqn} and \eqref{eq:adjEnergyEqn}. The computed fields are used to compute gradients and also as initial guess for the linear system at next optimization iteration. It turns out that it outperforms the previous strategy. \\
In this section, we show that model order reduction can be combined advantageously with the one-shot method to achieve even greater speedups. Similarly to the previous section, we consider different possible values for  the number of MOR basis vectors $r$. However, the MOR bases are constructed using solutions generated by the one-shot method instead of accurate solution as in previous section. This results in the following strategies:
\begin{itemize}
    \item MGCG\_1: Only one-shot method\citep{amir2024one},
    \item MOR\_2\_MGCG\_1: MOR with $r=2$ and $\romtau=5\cdot10^{-6}$. The one-shot is used as FOM solver when the MOR solution is not accurate enough,
    \item MOR\_3\_MGCG\_1: MOR with $r=3$ and $\romtau=5\cdot10^{-6}$.The one-shot is used as FOM solver when the MOR solution is not accurate enough,
    \item MOR\_4\_MGCG\_1: MOR with $r=4$ and $\romtau=5\cdot10^{-6}$. The one-shot is used as FOM solver when the MOR solution is not accurate enough.
\end{itemize}
\begin{figure}[h]
	\centering
	\begin{subfigure}[b]{0.3\textwidth}
		\centering
		\includegraphics[width=\textwidth]{Figures/3D_NL_fine/MGCG_full_iter/paraview_screenshot/eps.png}
		\caption{MGCG}
	\end{subfigure}     
    \hfill \\
	\begin{subfigure}[b]{0.24\textwidth}
		\centering
		\includegraphics[width=\textwidth]{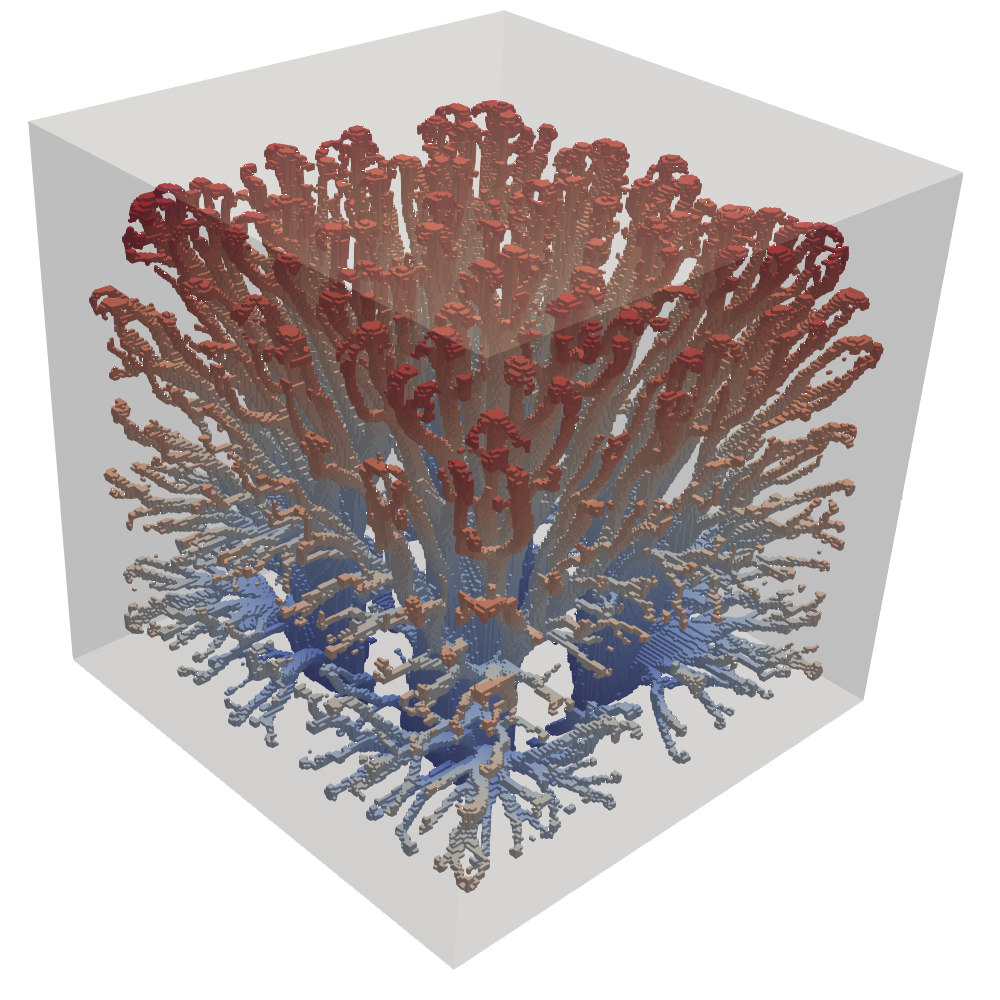}
		\caption{MGCG\_1}
	\end{subfigure}     
    \hfill
	\begin{subfigure}[b]{0.24\textwidth}
		\centering
		\includegraphics[width=\textwidth]{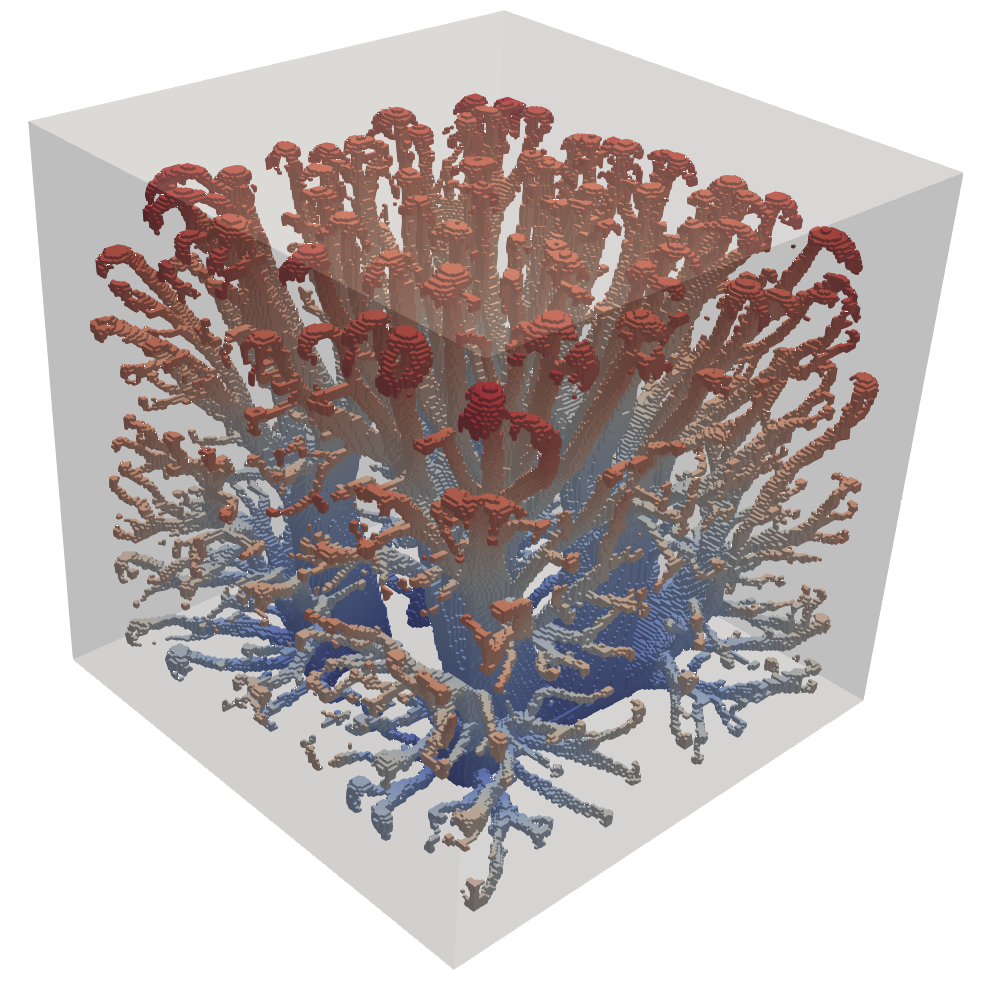}
		\caption{MOR\_2\_MGCG\_1}
	\end{subfigure}       
    \hfill 
	\begin{subfigure}[b]{0.24\textwidth}
		\centering
		\includegraphics[width=\textwidth]{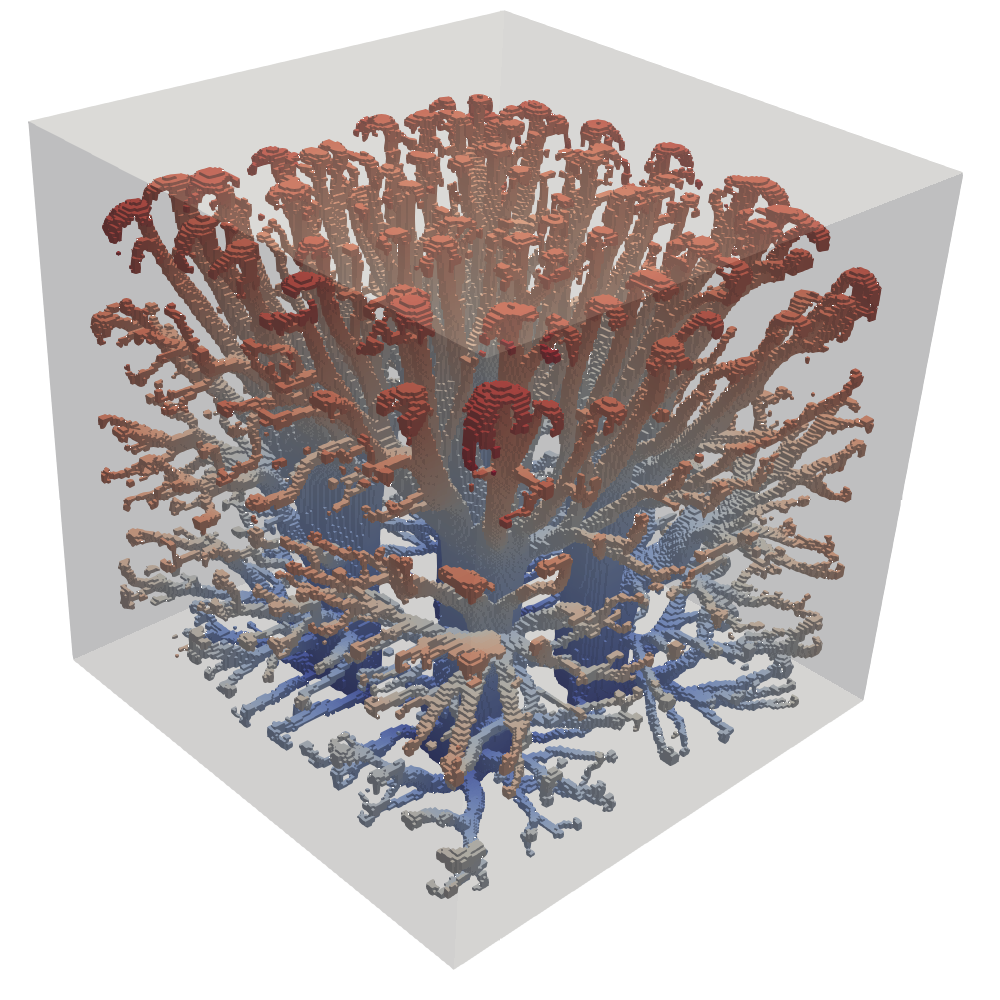}
		\caption{MOR\_3\_MGCG\_1}
	\end{subfigure}       
    \hfill
	\begin{subfigure}[b]{0.24\textwidth}
		\centering
		\includegraphics[width=\textwidth]{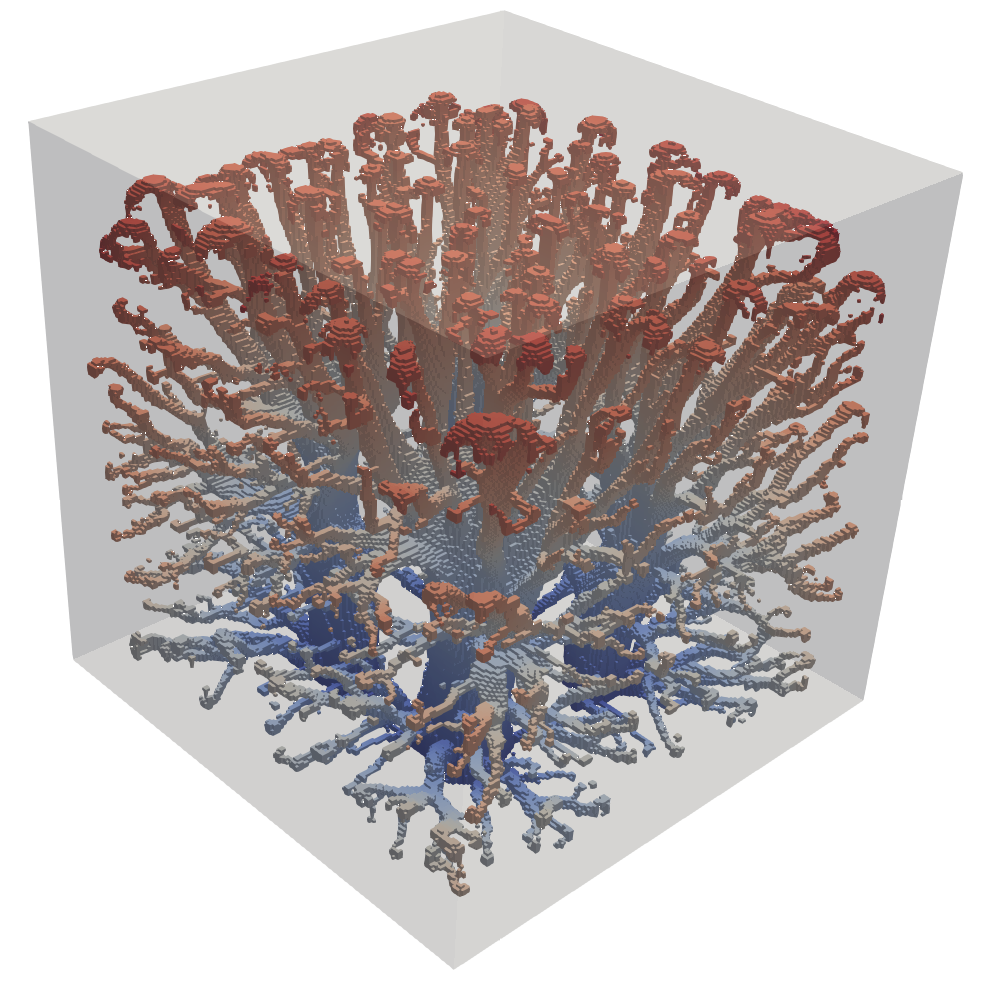}
	\caption{MOR\_4\_MGCG\_1}
	\end{subfigure}       
    \hfill
	\caption{Temperature field of the optimized design of the 3D thermal optimization using MOR. (a) Only accurate solver MGCG (b) one-shot method, (c)-(e) MOR with bases constructed using solutions generated by the one-shot method. A clip value of $0.5$ is used for the visualization}
	\label{fig:Inacc3D-2}
\end{figure}
\begin{figure}[h]
	\centering
    \begin{subfigure}[b]{0.45\textwidth}
		\centering
        \includegraphics[width=\textwidth]{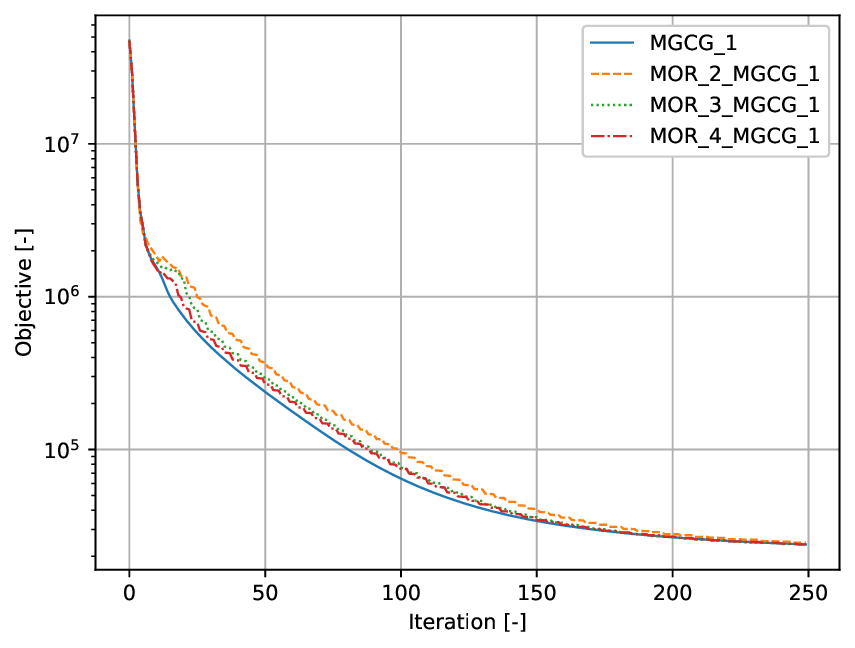}
        \caption{Objective value}
	\end{subfigure}
    \begin{subfigure}[b]{0.45\textwidth}
		\centering
        \includegraphics[width=\textwidth]{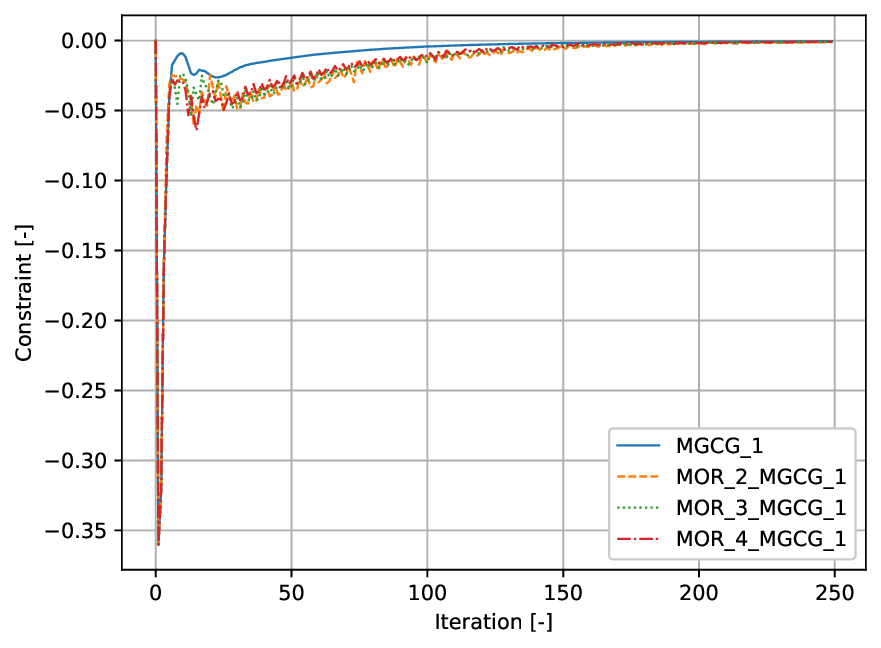}
        \caption{Constraint value}
	\end{subfigure}
    \caption{Objective value and constraint values along the optimization of the 3D thermal optimization using MOR in combination of the one-shot method.}
	\label{fig:MOR_inacc_obj_cons}
\end{figure}
The generated optimized designs for each of the inaccurate solver strategies are shown in figure \ref{fig:Inacc3D-2}, alongside the reference design optimized using MGCG until convergence. Due to the inaccuracies in the state field and the gradient, the optimized designs appear to be more asymmetrical than the reference MGCG design. However, the main tree-like structural features remain present.
Additionally, figure \ref{fig:MOR_inacc_obj_cons} plots the corresponding objective and constraint values for each of the inaccurate solver strategies. Here, there are visible differences between iteration $20$ and around $150$ where the one-shot method achieves lower objectives and constraint values closer to the goal. However, all strategies lead to similar objective values varying $2.41\%$ while satisfying the volume constraint.\\
To quantitatively assess the accuracy of each combination, we evaluate the relative error of the computed gradient compared to the accurate gradient (obtained with fully converged forward and adjoint equations) which is plotted in figure \ref{fig:3D_2_err}. The one-shot method achieves higher directional accuracy, evidenced by a lower angle deviation throughout all iterations. Moreover, its relative gradient error continuously decreases, suggesting progressively improving approximations.
The number of basis vectors does not seem to have a significant impact on the relative gradient error or the angle deviation since the curves are of the same order of magnitude. Nevertheless, these inaccuracies seem to have a limited impact on the final objective value, as seen in Table \ref{tab:InaccMGMOR}. All solver strategies achieve a final objective within $3\%$ of the one achieved by MGCG. All inaccurate workflows achieve a speedup with respect to the fully converged MGCG. Notably, the MOR workflows further reduce the total computational cost compared to the one-shot method, reaching a speedup of up to $1.57$. 
\begin{table}[ht]
  \centering
  \begin{tabular}{|l|c|c|c|c|c|}
    \hline
     & Objective & Walltime & Speedup & Forward    & Adjoint    \\ 
     & value     & (s)      & (CG)     & reductions & reductions \\ \hline
     MGCG\_1 & 23943.99 &  347.14  &  1& N/A & N/A \\ \hline
     MOR\_2\_MGCG\_1 & 24349.99 &  226.48  &  1.53& 181 & 173 \\ \hline
     MOR\_3\_MGCG\_1 & 23711.71 &  273.43  &  1.26& 183 & 173 \\ \hline
     MOR\_4\_MGCG\_1 & 23777.79 &  300.93  &  1.15& 183 & 172 \\ \hline 
  \end{tabular}
  \caption{Comparison of objective, timing, speedup and reduction counts}
  \label{tab:InaccMGMOR}
\end{table}
Table \ref{tab:InaccMGMOR} summarizes these performance outcomes, demonstrating the potential of MOR to enhance the optimization efficiency without significantly compromising the quality of the optimized design. 
\begin{table}[ht]
    \renewcommand{\arraystretch}{1.5}
    \centering
    \begin{tabular}{|c|c|c|}
    \hline
    \textbf{Solver name} & \textbf{Total time} & \textbf{Speedup}\\  \hline
    MGCG\_1 & 0.69 & 1\\ \hline
    MOR\_2 & 0.33 & 2.09\\ \hline
    MOR\_3 & 0.44 & 1.56 \\ \hline
    MOR\_4 & 0.53 & 1.30\\ \hline
    \end{tabular}
    \caption{Time of a single iteration of every solution strategy used for the 3D problem.}
    \label{tab:profile_fine}
\end{table}
The one-shot method relies heavily on the multigrid preconditioner to achieve satisfactory results with large computational gains. A single iteration of the one-shot method remains computationally more expensive than one MOR step, provided the basis size remains moderate. Therefore, we compare the computational efficiency of the MOR and one-shot methods by analyzing the cost of their core operations. For the one-shot method, these consist of the preconditioner-related operations and one conjugate gradient iteration. For MOR, these include creating and solving the reduced linear system, updating the basis, and performing related operations such as residual computation. Table \ref{tab:profile_fine} presents the time of one iteration of every solution strategy as well as the speedup with respect to the one one-shot iteration, calculated as $\nicefrac{t_{\text{one-shot}}}{t_{\text{MOR}}}$. Note that in the MOR strategies, the basis update step is executed only when the approximate solution fails to satisfy the residual tolerance. This conditional execution is accounted for in the average cost calculations. \\
The three MOR solvers are faster than the one-shot iteration, with a speedup of $2.09$, $1.56$ and $1.3$, respectively. As expected, the speedup decreases as the basis size increases. 
Note that the solver speedups shown in Table~\ref{tab:profile_fine} are higher than the effective global speedup shown in Table~\ref{tab:InaccMGMOR}. This discrepancy exists because Table~\ref{tab:profile_fine} strictly compares the isolated cost of a single MOR iteration against a single one-shot iteration and it does not account for the overhead of the fallback mechanism. Specifically, when the MOR approximation does not meet the residual threshold, the workflow incurs the computational cost of the failed MOR solver plus the cost of the one-shot method. Table~\ref{tab:InaccMGMOR} shows the total time spent solving the linear systems, either by the MOR solver or by the one-shot solver, thereby including the overhead of the fallback mechanism. 
It is worth noting that the one-shot method is a very high-performance method offering a speed-up of $10$ compared to the MGCG reference without the extra complexity added from the MOR operations. Nevertheless, an even greater speedup can be achieved when combining MOR with one-shot solutions.
\section{Conclusion}\label{s:conclusion}
In conclusion, this study investigated an approach to mitigate the high computational cost of density-based topology optimization by developing a framework that integrates model order reduction (MOR) to accelerate the state and adjoint solvers. Our approach is based on constructing two reduced models, one for the state and one for adjoint equations, to replace the computationally expensive, large linear systems. \\
Key contributions included a method for updating the reduced basis that is orthogonal and spans the same space as high-fidelity solutions of previous design iterations. Additionally, an alternative stopping criterion for the linear solvers and for the reduced models was investigated. This proved especially effective for the adjoint system.
To validate the performance, we presented a 3D heat sink design problem, where the proposed framework demonstrated significant performance gains. Two approaches to build the reduced bases were investigated: one consisting of converged state fields and  one consisting of one-shot solutions. In the former case, MOR achieved similar objective performance with high computational gains. In the latter case, when comparing to the one-shot method, computational gains were achieved with the same qualitative design and similar thermal performance. Furthermore, MOR offers a great potential for situations where the one shot method could not be effective such as more complex physics.
\section*{Declarations}
\textbf{Funding:} The work of Luis Cusicanqui was supported by  Vlaams Agentschap Innoveren en Ondernemen (VLAIO) through the Baekeland Research Project HBC.2022.0187. \\ \\
\textbf{Conflict of interest:}
On behalf of all authors, the corresponding author
states that there is no conflict of interest. \\ \\
\bibliographystyle{spbasic}      
\bibliography{bibliography}   
\end{document}

%% file: workflowDiagram.tex
\begin{tikzpicture}[node distance=0.8cm, scale=0.5]
    % Nodes
    \node (start) [startstop] {\scriptsize Start};
    \node (step1) [process, below of=start] {\scriptsize Initialize design field};
    \node (step2) [process, below of=step1, yshift=-0.3cm] {\scriptsize Solve energy equation \eqref{eq:energyEqn}};
    \node (step3) [process, below of=step2, yshift=-0.1cm] {\scriptsize Solve adjoint energy equation \eqref{eq:adjEnergyEqn}};
    \node (step4) [process, below of=step3, yshift=-0.1cm] {\scriptsize Compute gradients};
    \node (step5) [process, below of=step4, yshift=-0.1cm] {\scriptsize Update design variables};
    \node (step6) [decision, below of=step5, yshift=-0.3cm] {\tiny Converged?};
    \node (stop) [startstop, below of=step6, yshift=-0.3cm] {\scriptsize Optimum design};
    
    % Arrows
    \draw [arrow] (start) -- (step1);
    \draw [arrow] (step1) -- (step2);
    \draw [arrow] (step2) -- (step3);
    \draw [arrow] (step3) -- (step4);
    \draw [arrow] (step4) -- (step5);
    \draw [arrow] (step5) -- (step6);		
    \draw [arrow] (step6) -- node[anchor=east, xshift=0.5cm] {\tiny Yes} (stop);
    
    \node[fit=(step2) (step3), draw, dashed, inner sep=0.15cm](group){};
    \node at (group.north west) [right, yshift=0.2cm] {\tiny Solve FOM};
    
    \path (step1) -- (step2) coordinate[midway] (midpoint);
    
    \draw [arrow] (step6.east) node[yshift=-0.1cm, xshift=0.2cm] {\tiny No} -| ++(3,0) |- (midpoint);
\end{tikzpicture}	

%% file: mor_workflow_diagram.tex
\begin{tikzpicture}[node distance=0.8cm, scale=0.5]
	% Nodes
	\node (start) [startstop] {\scriptsize Start};
	\node (step1) [process, below of=start] {\scriptsize Initialize design field};
	\node (step2) [process, below of=step1, yshift=-0.3cm] {\scriptsize Solve energy equation \eqref{eq:energyEqn}};
	\node (step3) [process, below of=step2, yshift=-0.1cm] {\scriptsize Solve adjoint energy equation \eqref{eq:adjEnergyEqn}};
	\node (morstep1) [morprocess, below of=step3] {\scriptsize Update MOR};
	\node (step4) [process, below of=morstep1, yshift=-0.1cm] {\scriptsize Compute gradients};
	\node (step5) [process, below of=step4, yshift=-0.1cm] {\scriptsize Update design variables};
	\node (step6) [decision, below of=step5, yshift=-0.3cm] {\tiny Converged?};
	\node (stop) [startstop, below of=step6, yshift=-0.3cm] {\scriptsize Optimum design};
	\node (morstep2) [morprocess, right=of step6] {\scriptsize Solve MOR};
	\node (morstep3) [decision2, above=of morstep2, yshift=1.1cm] {\tiny Assess MOR};
	\path (step1) -- (step2) coordinate[midway] (mid12);
	\path (step4) -- (step5) coordinate[midway] (mid45);
	\path (morstep1) -- (step4) coordinate[midway] (midmor14);

	\draw [arrow] (start) -- (step1);
	\draw [arrow] (step1) -- (step2);
	\draw [arrow] (step2) -- (step3);
	\draw [arrow] (step3) -- (morstep1);
	\draw [arrow] (morstep1) -- (step4);
	\draw [arrow] (step4) -- (step5);
	\draw [arrow] (step5) -- (step6);	
	\draw [arrow] (morstep2) -- (morstep3);

	\draw [arrow] (step6) -- node[anchor=east, xshift=0.5cm] {\tiny Yes} (stop);
	
	\node[fit=(step2) (step3), draw, dashed, inner sep=0.15cm](group){};
	\node at (group.north west) [right, yshift=0.2cm] {\tiny Solve FOM};

	\draw [arrow] (step6.east) node[yshift=-0.1cm, xshift=0.2cm] {\tiny No} -- (morstep2.west);
	\draw [arrow] (morstep3.west) node[yshift=-0.1cm, xshift=-0.2cm] {\tiny Yes} -- (midmor14);
	\draw [arrow] (morstep3.north) node[yshift=0.1cm, xshift=-0.2cm] {\tiny No} -|++(0,2)|- (mid12);
\end{tikzpicture}	

%% file: figure_geometry.tex
\begin{tikzpicture}[scale=0.5]

\draw[thick] (0,0) rectangle (6,6);

\draw[ultra thick,blue] (2,6) -- (4,6) node[midway,above] {300K};

\draw[ultra thick,red, dashed] (0,6) -- (2,6) node[midway,right] {};
\draw[ultra thick,red, dashed] (4,6) -- (6,6) node[midway,right] {};
\draw[ultra thick,red,dashed] (6,6) -- (6,0) node[midway,right] {};
\draw[ultra thick,red,dashed] (6,0) -- (0,0) node[midway,below] {};
\draw[ultra thick,red,dashed] (0,0) -- (0,6) node[midway,below] {};

\begin{scope}
    \clip (0,0) -- (6,0) -- (6,-0.2) -- (0,-0.2) -- cycle; % Clip for the bottom side
    \foreach \x in {0,0.2,...,6} {
        \draw[thick,red] (\x,-0.2) -- ++(-0.2,0.2);
    }
\end{scope}
\draw[ultra thick,red,dashed] (6,0) -- (0,0) node[midway,below] {};
\begin{scope}
    \clip (6,0) -- (6,6) -- (6.2,6) -- (6.2,0) -- cycle; 
    \foreach \y in {0,0.2,...,6} {
        \draw[thick,red] (6.2,\y) -- ++(-0.2,0.2);
    }
\end{scope}
\begin{scope}
    \clip (-0.2,6) -- (-0.2,0) -- (0,0) -- (6,0) -- cycle; 
    \foreach \y in {0,0.2,...,6} {
        \draw[thick,red] (-0.2,\y) -- ++(0.2,0.2);
    }
\end{scope}
\draw[ultra thick,red,dashed] (6,6) -- (6,0) node[midway,right] {};

\end{tikzpicture}

%% file: 3D_geometry.tex
\begin{tikzpicture}[tdplot_main_coords, scale=1]  
  \coordinate (O)   at (0,0,0);
  \coordinate (X)   at (3,0,0);
  \coordinate (Y)   at (0,3,0);
  \coordinate (XY)  at (3,3,0);
  \coordinate (Z)   at (0,0,3);
  \coordinate (XZ)  at (3,0,3);
  \coordinate (YZ)  at (0,3,3);
  \coordinate (XYZ) at (3,3,3);
  
  \coordinate (0c)  at (1,1,0);
  \coordinate (Xc)  at (2,1,0);
  \coordinate (XYc) at (2,2,0);
  \coordinate (Yc)  at (1,2,0);

  \fill[red] (0c) -- (Xc) -- (XYc) -- (Yc) -- cycle;

  \node[black,anchor=center] at (1.5,1.5,0.02) {\tiny$\partial\Omega_d$};
  \node[black,anchor=center] at (2.5,2.5    ,0.02) {\tiny$\partial\Omega_n$};

  \draw[thick] (O) -- (X) -- (XY) -- (Y) -- cycle;      % bottom square
  \draw[thick] (Z) -- (XZ) -- (XYZ) -- (YZ) -- cycle;   % top square
  \draw[thick] (O) -- (Z);
  \draw[thick] (X) -- (XZ);
  \draw[thick] (XY) -- (XYZ);
  \draw[thick] (Y) -- (YZ);

  \coordinate (AxisOrigin) at (3.5, 3.5, 1.5);
  \draw[->, thick] (AxisOrigin) -- ++(0.4, 0, 0) 
    node[anchor=north west]{$x$};
  \draw[->, thick] (AxisOrigin) -- ++(0, 0.4, 0) 
    node[anchor=south]{$y$};
  \draw[->, thick] (AxisOrigin) -- ++(0, 0, 0.4) 
    node[anchor=south]{$z$};
\end{tikzpicture}